\documentclass[hidelinks,onefignum,onetabnum]{siamart251216}

\usepackage{lipsum}
\usepackage{hyperref}
\usepackage{amsfonts}
\usepackage{graphicx}
\usepackage{epstopdf}
\usepackage{bm}
\usepackage{amsmath}
\usepackage{amssymb}
\usepackage{mathrsfs}
\usepackage{enumerate}
\usepackage{placeins}
\usepackage{color}
\usepackage{subcaption}
\usepackage{algorithmic}
\usepackage[normalem]{ulem} 

\ifpdf
\DeclareGraphicsExtensions{.eps,.pdf,.png,.jpg}
\else
  \DeclareGraphicsExtensions{.eps}
\fi

\newsiamremark{rem}{Remark}
\newsiamremark{hypothesis}{Hypothesis}
\crefname{hypothesis}{Hypothesis}{Hypotheses}
\newsiamthm{claim}{Claim}
\newsiamthm{assm}{Assumption}
\newsiamthm{cdn}{Condition}
\newsiamthm{prop}{Proposition}
\newsiamremark{fact}{Fact}
\crefname{fact}{Fact}{Facts}

\headers{}{T. Torabi, C. H. Ho, C. Ortner, and R. N. Qiu}

\title{Recovering from Averages: An Inverse Problem Perspective on Machine Learning Interatomic Potentials}

\author{
Tina Torabi$^{\dagger}$\thanks{Department of Mathematics,
University of British Columbia,
Vancouver, BC V6T 1Z2, Canada.
(\email{torabit@math.ubc.ca}, corresponding author)}
\and
Cheuk Hin Ho\thanks{Department of Mathematics,
University of British Columbia,
Vancouver, BC V6T 1Z2, Canada.}
\and
Christoph Ortner\footnotemark[2]
\and
Ruo Ning Qiu\thanks{Department of Mathematics,
University of California San Diego,
La Jolla, CA, USA; and
Department of Mathematics and Statistics,
San Diego State University,
San Diego, CA, USA.}
}

\begin{document}

\maketitle

\begin{abstract}
Machine learning interatomic potentials (MLIPs) approximate the potential
energy surface of atomistic systems from quantum-mechanical reference data
and enable molecular simulation at near-ab-initio accuracy at a fraction of
the cost. Since they are fitted to total energies and forces, it is not
obvious what is estimated: the total energy or the many-body potentials from which it is assembled.
Mathematically, this can be formulated as a question whether a 
function can be recovered from averages of its values over random clusters. 
We study it for a linear
model of fixed correlation order in simplified statistical setting. We
show that, in the population limit, the potential is identifiable in a
natural weighted Sobolev space, and derive finite-sample recovery
guarantees.
Numerical experiments for pair, three-body, and four-body models support
the theory.
\end{abstract}



\section{Introduction}
Machine learning interatomic potentials (MLIPs)~\cite{Behler2007,
Drautz2019} approximate the potential energy surface of an atomistic
system by a model fitted to quantum-mechanical reference data, typically
obtained from density functional theory (DFT)~\cite{Hohenberg1964}. They
reach near-DFT accuracy at a small fraction of the cost and have become a
standard tool in molecular simulation. The state of the art is dominated by
equivariant graph neural networks (EGNNs)~\cite{Batzner_2022,
NEURIPS2022_4a36c3c5}, but even linear models such as the moment tensor
potential~\cite{doi:10.1137/15M1054183}, SNAP~\cite{THOMPSON2015316}, and the
atomic cluster expansion~\cite{Drautz2019, DUSSON2022110946} remain
competitive~\cite{LACE} for many applications. 

These linear MLIPs express the total energy of a configuration $\mathbf R$ as a many-body expansion (some but not all EGNNs do so {\it implicitly})
\[
    E(\mathbf R)
    =
    \sum_{\nu=1}^{\nu_{\max}}
    \sum_{c\in\widetilde{\mathcal P}_\nu(\mathbf R)}
    \widetilde V^{(\nu)}(c),
\]
where $\widetilde{\mathcal P}_\nu(\mathbf R)$ is the set of clusters of $\nu+1$
atoms within a cutoff radius and $\widetilde V^{(\nu)}$ is a $(\nu+1)$-body
potential. The training data consist of total energies $E(\mathbf R)$ and forces
$-\nabla E(\mathbf R)$ for a set of configurations. Each observation is
therefore a sum of many evaluations of the unknown potentials $\widetilde
V^{(\nu)}$. We refer to the problem of estimating the $\widetilde V^{(\nu)}$
from such data as \emph{recovery from averages}. 

This raises the question which object a fit to energies and forces identifies.
The total energy is observed directly, and predictive accuracy for $E$ is the
standard measure of success~\cite{Kov_cs_2023, LACE, 2025-sparse-ace-1}. The
potentials $\widetilde V^{(\nu)}$ are not observed. A fit identifies them only
if the map from the potentials to the observations is injective. This
distinction matters in practice: in simulation, the behaviour of the model,
especially out-of-distribution, is determined by the potentials and their
derivatives, not by its accuracy on the training labels. The questions we
address are whether, and in which norm, the potentials $\widetilde
V^{(\nu)}$ can be recovered from energy and force observations, and
whether the averaging over many clusters in each observation comes at a
price of requiring more data when each observation aggregates
more clusters?

We study this in a simplified setting: a linear model of a single
correlation order $\nu$, and training configurations whose clusters are
independent samples from a cluster distribution $\rho$. Energy observations
are then averages of $\widetilde V^{(\nu)}$ over the $K$ clusters of a
configuration, and force observations are averages of its derivatives. Our results are as follows: In the
population limit, the least-squares loss for energy observations is an
$L^2(\rho)$-type norm of the error.
Adding force observations adds control of the derivatives; the loss
then becomes a weighted Sobolev-type norm. 
These results demonstrate that $\widetilde V^{(\nu)}$ is indeed identifiable
from aggregate observations. 
Finite-sample recovery guarantees follow from the
general framework of Adcock et al.~\cite{pmlr-v235-adcock24a}. The
worst-case guarantees grow with the number of clusters per observation,
but the expected observation size does not, and numerical experiments
for pair, three-body, and four-body potentials confirm the latter. 

\paragraph{Related work}
Our analysis is a direct application of the theory of least-squares
approximation from random samples~\cite{cohen2018, cohen2017optimal,
pmlr-v235-adcock24a}, which admits general linear observation operators and
several observation types. The novelty of our work lies in the observation
operators: each sample is an average over many clusters, and the resulting
population geometry and variation differ from those of point evaluations. The
results apply whenever a global quantity is learned from a sum of local
contributions. The machine learning literature on learning from aggregate
observations~\cite{3495724.3496394} treats a related problem: Zhang et al.\
characterize when maximum-likelihood estimators recover individual-level
parameters from aggregate supervision, up to the equivalence induced by the
aggregation. We instead ask which norm governs stable recovery and how the
required number of configurations depends on the number of clusters per
observation and on the model dimension. In the MLIP literature the questions of
identifiability and stability have received little attention; the emphasis is on
predictive accuracy~\cite{Kov_cs_2023, LACE, 2025-sparse-ace-1}.

\paragraph{Scope and limitations}
The analysis is restricted to a single correlation order, linear
approximation spaces, and independent clusters with a fixed, or bounded
random, number of clusters per observation. Realistic MLIPs combine
several body orders, and training data generated from molecular dynamics
trajectories are correlated. The single-order i.i.d.\ setting isolates the
recovery-from-averages structure and allows the induced geometry to be
characterized exactly and with elementary methods.

\paragraph{Outline}
\Cref{sec:invsform} formulates the problem: many-body potentials, energy and
force observations, the statistical model, and the least-squares estimator.
\Cref{sec:theory} identifies the population geometry and the Gram matrix, and
establishes finite-sample recovery guarantees. \Cref{sec:numerics} presents the
numerical experiments. The appendices contain the detailed derivation, in particular proofs of the variation bounds.

\section{Background and Problem Formulation}\label{sec:invsform}
\subsection{Many-body interatomic potentials}
An atomic configuration of $N_{\rm at}$ atoms is a set
$\mathbf R := \{\bm r_1,\dots,\bm r_{N_{\rm at}}\}\subset\mathbb R^d$
of atomic positions. An interatomic potential is a map
$\mathbf R\mapsto E(\mathbf R)$ that assigns a potential energy to each
configuration; it is used, for example, in molecular dynamics and geometry
optimization. (In practice, atoms are annotated by their chemical elements, 
which we ignore for the sake of notational simplicity.)

For $\nu\geq1$, a \emph{$\nu$-cluster} of $\mathbf R$ is a tuple
$c=(\bm r_{j_0},\dots,\bm r_{j_\nu})$, $j_0<\dots<j_\nu$, satisfying a
locality criterion, for example that all pairwise distances within $c$ are
at most a cutoff radius $r_{\rm cut}$. We denote the set of $\nu$-clusters
of $\mathbf R$ by $\widetilde{\mathcal P}_\nu(\mathbf R)$. Most
interatomic potentials assume that the total energy is a sum of
contributions of increasing correlation order,
\begin{equation}\label{eq:many_body_expansion}
    E(\mathbf R)
    =
    \sum_{\nu=1}^{\nu_{\max}}
    \sum_{c\in\widetilde{\mathcal P}_\nu(\mathbf R)}
    \widetilde V^{(\nu)}(c),
\end{equation}
where the $(\nu+1)$-body potential
$\widetilde V^{(\nu)}:(\mathbb R^d)^{\nu+1}\to\mathbb R$ is invariant under
the Euclidean group $E(d)$ and under permutations of its arguments.

By translation invariance we may place one atom of each cluster at the
origin. For $c\in\widetilde{\mathcal P}_\nu(\mathbf R)$ we define the
relative coordinates $\bm z^{(c)}\in(\mathbb R^d)^\nu$ and the reduced
potential $V^{(\nu)}:(\mathbb R^d)^\nu\to\mathbb R$ by
\[
    \bm z^{(c)}
    :=
    (\bm r_{j_1}-\bm r_{j_0},\dots,\bm r_{j_\nu}-\bm r_{j_0}),
    \qquad
    V^{(\nu)}(\bm z^{(c)})
    :=
    \widetilde V^{(\nu)}(c).
\]
By permutation invariance, any atom of the cluster may serve as the origin.
The reduced potential remains invariant under the orthogonal group $O(d)$; in
particular, $V^{(\nu)}(-\bm z)=V^{(\nu)}(\bm z)$. Collecting the relative
coordinates in the multiset
$\mathcal Z_\nu(\mathbf R):=\{\bm z^{(c)}:c\in\widetilde{\mathcal P}_\nu(\mathbf R)\}$,
the expansion \eqref{eq:many_body_expansion} becomes
\[
    E(\mathbf R)
    =
    \sum_{\nu=1}^{\nu_{\max}}
    \sum_{\bm z\in\mathcal Z_\nu(\mathbf R)}
    V^{(\nu)}(\bm z).
\]

Our aim is to estimate the potentials $V^{(\nu)}$ from observations about the 
total energy $E$. We focus on linear parametrizations: models such
as the Moment Tensor Potential (MTP) and the Atomic
Cluster Expansion (ACE)~\cite{Drautz2019, DUSSON2022110946,
doi:10.1137/15M1054183} parametrize each $V^{(\nu)}$ as
\[
    V^{(\nu)}_\theta(\bm z)
    =
    \sum_{\varphi\in B_\nu}\theta_\varphi\,\varphi(\bm z),
\]
where $B_\nu$ is a set of $O(d)$- and permutation-invariant basis functions
$\varphi:(\mathbb R^d)^\nu\to\mathbb R$. The resulting total energy
\begin{equation}\label{eq:many_body_exp}
    E_\theta(\mathbf R)
    =
    \sum_{\nu=1}^{\nu_{\max}}
    \sum_{\bm z\in\mathcal Z_\nu(\mathbf R)}
    V^{(\nu)}_\theta(\bm z)
\end{equation}
is linear in the parameters $\theta$.

\subsection{Parameter estimation}
Let $\{\mathbf R_\alpha\}_{\alpha=1}^{N}$ be a training set of $N$ atomic
configurations, and suppose that for each $\mathbf R_\alpha$ reference
calculations provide a total energy $E_\alpha^{\rm ref}$ and forces
$\bm F_{\alpha,i}^{\rm ref}$, $i=1,\ldots,N_{\rm at}^{(\alpha)}$. The model
parameters are commonly estimated by minimizing the weighted least-squares
loss
\begin{equation}\label{eq:energy_force_loss} 
    \widehat{\theta}
    \in
    \operatorname*{arg\,min}_{\theta}
    \frac{1}{N}
    \sum_{\alpha=1}^{N}
    \left[
        w_{\rm E}
        \left|
            E_\theta(\mathbf R_\alpha)-E_\alpha^{\rm ref}
        \right|^2
        +
        w_{\rm F}
        \sum_{i=1}^{N_{\rm at}^{(\alpha)}}
        \left\|
            -\nabla_{\bm r_i}E_\theta(\mathbf R_\alpha)
            -
            \bm F_{\alpha,i}^{\rm ref}
        \right\|_2^2
    \right],
\end{equation}
where $w_{\rm E},w_{\rm F}\geq 0$ weight the energy and force
contributions.

The inverse problem underlying \eqref{eq:energy_force_loss} is to determine the
potentials $V^{(\nu)}$, or equivalently the parameters $\theta$, from total
energy and force observations. In the sense of Hadamard, this
problem is well posed if a solution exists, is unique, and depends continuously
on the data. Existence is not a concern for us; uniqueness and stability are:  
First-principles calculations provide only total energies and
forces, each of which aggregates many cluster contributions, and it is
conceivable that different potentials $V^{(\nu)}$ produce identical or nearly
identical observations. This raises the central question of this work: can the
potentials be recovered from such aggregated observations and, if so, is the
recovery stable in a natural norm?

\begin{rem}[Non-identifiability across body orders]
\label{rem:mixed_orders}
Uniqueness can fail when several body orders are fitted simultaneously.
Consider $\nu_{\max}=2$ and configurations of $n\geq3$ atoms in which all
pairs and triples lie within the cutoff, so that
$\#\mathcal Z_1(\mathbf R)=\binom n2$ and
$\#\mathcal Z_2(\mathbf R)=\binom n3$. For any constant $c$, replacing
$V^{(1)}$ by $V^{(1)}+c$ and $V^{(2)}$ by $V^{(2)}-3c/(n-2)$ changes the
total energy by
$c\binom n2-\frac{3c}{n-2}\binom n3=0$
and leaves all forces unchanged. The two models produce identical
observations on every such configuration, although their components
differ. For configurations with varying numbers of clusters the
cancellation is no longer exact, but the same direction in parameter
space remains nearly invisible to the data. We therefore consider a single
correlation order throughout.
\end{rem}

\subsection{Simplified statistical model}

In view of \cref{rem:mixed_orders}, we fix a single correlation order
$\nu$ and set $V^{(\nu')}\equiv0$ for $\nu'\neq\nu$. We write $V:=V^{(\nu)}$,
$V_\theta:=V^{(\nu)}_\theta$, and $\mathcal Z:=\mathcal Z_\nu$.

Let $\mathcal Z$ be the collection of cluster coordinates of a
configuration and, for an atom $i$, let $\mathcal Z'$ be the collection of
clusters containing $i$, with $\bm r_i$ placed at the origin. The total
energy and the force on atom $i$ are then
\[
    E=\sum_{\bm z\in\mathcal Z}V(\bm z),
    \qquad
    \bm F_i=\sum_{\bm z\in\mathcal Z'}DV(\bm z),
    \qquad
    DV(\bm z):=\sum_{\ell=1}^{\nu}\nabla_{\bm z_\ell}V(\bm z),
\]
where $\bm z=(\bm z_1,\ldots,\bm z_\nu)\in(\mathbb R^d)^\nu$; the force is
derived in Appendix~\ref{sec:force_obsv}. To simplify the presentation, we
assume that the numbers of clusters are fixed.

\begin{assm}[Fixed numbers of clusters]\label{assm:fixedK}
There exist $K,K'\in\mathbb N$ such that $\#\mathcal Z=K$ and
$\#\mathcal Z'=K'$ for all configurations and atoms.
\end{assm}

The extension to a bounded random number of clusters is straightforward
by conditioning on $K$; $K$ is then replaced by $\mathbb E[K]$ in the
results of Section~\ref{sec:theory} and by $K_{\max}$ in the variation
bounds. Typically $K'\ll K$, since $\mathcal Z'$ contains only the
clusters involving atom $i$. We define the normalized energy and force
observation operators
\begin{equation}\label{eq:norm_obs}
    \mathcal O_{\rm E}(V;\mathcal Z)
    :=
    \frac{1}{\sqrt K}\sum_{\bm z\in\mathcal Z}V(\bm z),
    \qquad
    \mathcal O_{\rm F}(V;\mathcal Z')
    :=
    \frac{1}{\sqrt{K'}}\sum_{\bm z\in\mathcal Z'}DV(\bm z).
\end{equation}
The scaling normalizes the variance of the observations; see
Lemma~\ref{lem:population_observation_geometry}.

Assuming noise-free reference data generated by $V$, and up to a
rescaling of the weights, the loss \eqref{eq:energy_force_loss} can be
written in terms of the observation operators. Let
$(\mathcal Z_\alpha)_{\alpha=1}^{N}$ and $(\mathcal Z'_\beta)_{\beta=1}^{N'}$
denote the cluster collections of the energy and force observations,
respectively. By
linearity of the observation operators, the empirical loss is 

\begin{equation}\label{eq:abstract_loss}
\resizebox{0.95\linewidth}{!}{$\displaystyle
    \mathcal L_{N,N'}(V_\theta)
    =     \frac{w_{\rm E}}{N}
    \sum_{\alpha=1}^{N}
    \left|
        \mathcal O_{\rm E}(V_\theta;\mathcal Z_\alpha)
        -
        \mathcal O_{\rm E}(V;\mathcal Z_\alpha)
    \right|^2+
    \frac{w_{\rm F}}{N'}
    \sum_{\beta=1}^{N'}
    \left\|
        \mathcal O_{\rm F}(V_\theta;\mathcal Z'_\beta)
        -
        \mathcal O_{\rm F}(V;\mathcal Z'_\beta)
    \right\|_2^2
$}
\end{equation}

The least-squares estimator is
 $\widehat{\theta}
    \in
    \operatorname*{arg\,min}_{\theta}
    \mathcal{L}_{N,N'}(V_\theta).$

\begin{rem}[Relative weighting of energy and force data]
Suppose that $N'/N\approx\gamma>0$. Then
\[
    \frac{\bar w_{\rm E}}{\bar w_{\rm F}}
    \approx \gamma\frac{w_{\rm E}}{w_{\rm F}}
    =
    \frac{N'}{N}\frac{w_{\rm E}}{w_{\rm F}}.
\]
Thus, since $N\ll N'$, this normalization assigns
a larger effective weight to each energy observation. This compensates for
the smaller number of energy observations and is consistent with the larger
per-observation energy weights commonly used in practice.
\end{rem}
    
In this form, the inverse problem is a \emph{recovery-from-averages}
problem: each observation is a scaled sum of evaluations of $V$, or of
$DV$, over many clusters. We adopt a statistical model in which the
cluster coordinates are random. 

\begin{assm}[Statistical model]\label{assm:indep}
There are probability measures $\rho$ and $\rho'$ on $(\mathbb R^d)^\nu$
such that
\begin{enumerate}
    \item \textbf{Cluster independence:} the cluster coordinates of an
    energy observation are i.i.d.\ with law $\rho$, and those of a force
    observation are i.i.d.\ with law $\rho'$; that is,
    $\mathcal Z\sim\rho^{\otimes K}$ and $\mathcal Z'\sim\rho'^{\otimes K'}$;
    \item \textbf{Inversion symmetry:} $\rho'$ is invariant under
    $\bm z\mapsto-\bm z$.
\end{enumerate}
The energy samples $\mathcal Z_\alpha$ are i.i.d.\ draws from
$\rho^{\otimes K}$, and the force samples $\mathcal Z'_\beta$ are i.i.d.\
draws from $\rho'^{\otimes K'}$, independent of the $\mathcal Z_\alpha$.
\end{assm}

\begin{rem}[Configuration-level sampling]
Assumption~\ref{assm:indep} is an idealization. Clusters within a
configuration share atoms and are therefore correlated. Moreover, in
\eqref{eq:energy_force_loss} the forces are evaluated on the same
configurations as the energies, so that $N'=\sum_\alpha N_{\rm at}^{(\alpha)}$
and the $\mathcal Z'_\beta$ are not independent of the $\mathcal Z_\alpha$.
The measures $\rho$ and $\rho'$ differ in general, since $\mathcal Z'$
contains only clusters that involve a given atom. Inversion symmetry holds,
for example, whenever $\rho'$ is invariant under $O(d)$.
A more realistic statistical model would first sample the atomic positions and
then construct the associated clusters. For example, the atoms could be
generated by a Poisson point process
~\cite{642efb3b46d5467caa8b384a27968448}. In this formulation, the number of
clusters is random, and clusters sharing atoms are dependent. Extending the results below to such configuration-level models is left to
future work.
\end{rem}

\subsection{Linear approximation and least squares}\label{sec:LSQframework}

We enumerate the basis as $B_\nu=\{\varphi_1,\dots,\varphi_M\}$ and write
$U_M:=\operatorname{span}B_\nu$, so that
\[
    V_\theta=\sum_{m=1}^M\theta_m\varphi_m,
    \qquad
    \theta=(\theta_1,\dots,\theta_M)^\top\in\mathbb R^M.
\]
For $m=1,\dots,M$, we define the energy and force features
\begin{equation}\label{eq:features}
    \begin{aligned}
    \Phi_m^{\rm E}(\mathcal Z)
    &:=
    \mathcal O_{\rm E}(\varphi_m;\mathcal Z)
    =
    \frac{1}{\sqrt K}\sum_{\bm z\in\mathcal Z}\varphi_m(\bm z),
    \qquad \text{and} \\ 
    \Phi_m^{\rm F}(\mathcal Z')
    &:=
    \mathcal O_{\rm F}(\varphi_m;\mathcal Z')
    =
    \frac{1}{\sqrt{K'}}\sum_{\bm z\in\mathcal Z'}D\varphi_m(\bm z).
    \end{aligned}
\end{equation}
By linearity,
$\mathcal O_{\rm E}(V_\theta;\mathcal Z)=\sum_{m}\theta_m\Phi_m^{\rm E}(\mathcal Z)$
and
$\mathcal O_{\rm F}(V_\theta;\mathcal Z')=\sum_{m}\theta_m\Phi_m^{\rm F}(\mathcal Z')$.

The least-squares problem \eqref{eq:abstract_loss} is therefore
equivalent to the normal equations $G\theta=f$, where
\begin{equation}\label{eq:empirical_gram}
\begin{aligned}
    G_{mn}
    &=
    \frac{w_{\rm E}}{N}\sum_{\alpha=1}^{N}
    \Phi_m^{\rm E}(\mathcal Z_\alpha)\Phi_n^{\rm E}(\mathcal Z_\alpha)
    +
    \frac{w_{\rm F}}{N'}\sum_{\beta=1}^{N'}
    \Phi_m^{\rm F}(\mathcal Z'_\beta)\cdot\Phi_n^{\rm F}(\mathcal Z'_\beta),
    \\
    f_m
    &=
    \frac{w_{\rm E}}{N}\sum_{\alpha=1}^{N}
    \Phi_m^{\rm E}(\mathcal Z_\alpha)\,\mathcal O_{\rm E}(V;\mathcal Z_\alpha)
    +
    \frac{w_{\rm F}}{N'}\sum_{\beta=1}^{N'}
    \Phi_m^{\rm F}(\mathcal Z'_\beta)\cdot\mathcal O_{\rm F}(V;\mathcal Z'_\beta).
\end{aligned}
\end{equation}
The parameters are uniquely determined if and only if $G$ is positive
definite, and the stability of the recovery is governed by the spectrum
of $G$. In Section~\ref{sec:theory} we characterize the population Gram
matrix $\mathbb E[G]$ and the norm it induces on $U_M$, and bound the
numbers of samples $N$ and $N'$ required for $G$ to inherit its
properties.

\section{Theoretical results}\label{sec:theory}

We first study the infinite-data limit of the least-squares
problem~\eqref{eq:abstract_loss}, which identifies the norms induced by
energy and force observations. We then restrict to the approximation
space $U_M$, where these norms are represented by the population Gram
matrix, and finally bound the number of samples required for the
empirical Gram matrix to inherit its properties.

\subsection{Population geometry}
Taking the expectation of the loss~\eqref{eq:abstract_loss} under
Assumption~\ref{assm:indep} gives the population loss
\begin{equation}\label{eq:population_loss}
    \mathcal L_\infty(V_\theta)
    :=
    w_{\rm E}\,
    \mathbb E_{\mathcal Z\sim\rho^{\otimes K}}
    \left|\mathcal O_{\rm E}(V_\theta-V;\mathcal Z)\right|^2
    +
    w_{\rm F}\,
    \mathbb E_{\mathcal Z'\sim\rho'^{\otimes K'}}
    \left\|\mathcal O_{\rm F}(V_\theta-V;\mathcal Z')\right\|_2^2,
\end{equation}
which is the limit of the empirical loss as $N,N'\to\infty$. The following
lemma computes the two expectations for a generic function $u$; in our
application, $u=V_\theta-V$.

\begin{lemma}[Geometry of the population observations]
\label{lem:population_observation_geometry}
Let Assumption~\ref{assm:indep} hold, and let $u:(\mathbb R^d)^\nu\to\mathbb R$
be $O(d)$-invariant with $u\in L^2(\rho)$ and $Du\in L^2(\rho';\mathbb R^d)$.
Then
\begin{align}
    \mathbb E_{\mathcal Z\sim\rho^{\otimes K}}
    \left|\mathcal O_{\rm E}(u;\mathcal Z)\right|^2
    &=
    \|u\|_{L^2(\rho)}^2
    +(K-1)\left|\int u\,d\rho\right|^2,
    \label{eq:expected_energy_observation}
    \\
    \mathbb E_{\mathcal Z'\sim\rho'^{\otimes K'}}
    \left\|\mathcal O_{\rm F}(u;\mathcal Z')\right\|_2^2
    &=
    \|Du\|_{L^2(\rho';\mathbb R^d)}^2.
    \label{eq:expected_force_observation}
\end{align}
\end{lemma}

\begin{proof}
Write $\mathcal Z=\{\bm z^{(1)},\ldots,\bm z^{(K)}\}$, counted with
multiplicity. Expanding the square and separating the $K$ diagonal from
the $K(K-1)$ off-diagonal terms,
\[
    \mathbb E\left|\mathcal O_{\rm E}(u;\mathcal Z)\right|^2
    =
    \frac1K\sum_{k=1}^K\mathbb E\,|u(\bm z^{(k)})|^2
    +\frac1K\sum_{k\neq j}\mathbb E\big[u(\bm z^{(k)})u(\bm z^{(j)})\big].
\]
Since the $\bm z^{(k)}$ are i.i.d.\ with law $\rho$, each diagonal term
equals $\|u\|_{L^2(\rho)}^2$ and each off-diagonal term equals
$\big(\int u\,d\rho\big)^2$, which gives
\eqref{eq:expected_energy_observation}.

The same expansion for $\mathcal Z'=\{\bm z^{(1)},\ldots,\bm z^{(K')}\}$
gives
\[
    \mathbb E\left\|\mathcal O_{\rm F}(u;\mathcal Z')\right\|_2^2
    =
    \|Du\|_{L^2(\rho';\mathbb R^d)}^2
    +(K'-1)\left\|\int Du\,d\rho'\right\|_2^2.
\]
Since $u$ is $O(d)$-invariant, $u(-\bm z)=u(\bm z)$, hence $Du$ is odd.
By assumption, $\rho'$ is invariant under $\bm z\mapsto-\bm z$, hence the change of
variables $\bm z\mapsto-\bm z$ gives
$\int Du\,d\rho'=-\int Du\,d\rho'$, so the mean term vanishes and
\eqref{eq:expected_force_observation} follows.
\end{proof}

Motivated by Lemma~\ref{lem:population_observation_geometry}, we define
the energy-induced and energy-force-induced norms
\[
    \|u\|_{\rm E}^2
    :=
    \|u\|_{L^2(\rho)}^2+(K-1)\left|\int u\,d\rho\right|^2,
    \qquad
    \|u\|_{\rm EF}^2
    :=
    w_{\rm E}\|u\|_{\rm E}^2+w_{\rm F}\|Du\|_{L^2(\rho';\mathbb R^d)}^2,
\]
so that $\mathcal L_\infty(V_\theta)=\|V_\theta-V\|_{\rm EF}^2$. Since
$|\int u\,d\rho|^2\leq\|u\|_{L^2(\rho)}^2$, the energy norm is equivalent
to the $L^2(\rho)$ norm,
$\|u\|_{L^2(\rho)}\leq\|u\|_{\rm E}\leq\sqrt K\,\|u\|_{L^2(\rho)}$,
and differs from it only through the weight on the mean. Energy
observations therefore control the error in $L^2(\rho)$, with the mean
weighted by $K$; force observations additionally control its derivative
in $L^2(\rho')$.

\begin{rem}[Origin-exchange symmetry and the $H^1$ norm] 
\label{rem:origin_exchange}
The force norm $\|Du\|_{L^2(\rho';\mathbb R^d)}$ appears to control only the
combination $Du=\sum_\ell\nabla_{\bm z_\ell}u$, which is the negative force on
the origin atom of the cluster. The full gradient $\nabla u=(\nabla_{\bm
z_\ell}u)_{\ell=1}^d$ collects the force contributions on the remaining $\nu$
atoms. Let $T_\ell$ denote the change of origin from atom $0$ to atom $\ell$,
\[
    T_\ell\bm z
    :=
    (\bm z_1-\bm z_\ell,\ldots,\bm z_{\ell-1}-\bm z_\ell,
     -\bm z_\ell,\bm z_{\ell+1}-\bm z_\ell,\ldots,\bm z_\nu-\bm z_\ell).
\]
By permutation invariance of $\widetilde V$, $\nabla_{\bm z_\ell}u(\bm
z)=-Du(T_\ell\bm z)$. If $\rho'$
is invariant under each $T_\ell$, i.e.\ that the atoms of a force cluster are
statistically indistinguishable, then 
\[
    \|\nabla u\|_{L^2(\rho';(\mathbb R^d)^\nu)}^2
    =
    \nu\,\|Du\|_{L^2(\rho';\mathbb R^d)}^2;
\]
that is, $\|D \cdot \|_{L^2(\rho')}$ becomes an $H^1$ semi-norm.
\end{rem}

\subsection{Population Gram matrix on the approximation space}
\label{sec:population_gram}

Lemma~\ref{lem:population_observation_geometry} does not depend on the
approximation space. We now restrict to $U_M$ and impose a normalization
of the basis. This has two consequences: the best-approximation error in
the induced norms reduces to the $L^2(\rho)$ and Sobolev-type errors, and
the induced norms are represented exactly by the population Gram matrix.

\begin{assm}[Basis normalization]
\label{assm:basis_normalization}
The first basis function is constant, $\varphi_1\equiv 1/\sqrt K$, and
$\varphi_2,\ldots,\varphi_M$ are orthonormal in $L^2(\rho)$ with zero
mean,
\[
    \int_{(\mathbb R^d)^\nu}\varphi_m(\bm z)\,d\rho(\bm z)=0,
    \qquad m=2,\ldots,M.
\]
\end{assm}

The energy norm differs from the $L^2(\rho)$ norm only in the weight on
the mean,
\[
    \|u\|_{\rm E}^2
    =
    \|u-\bar u\|_{L^2(\rho)}^2+K|\bar u|^2,
    \qquad
    \bar u:=\int u\,d\rho.
\]
Since $U_M$ contains the constants, the mean of the residual can be matched
exactly.

\begin{corollary}[Best-approximation error]
\label{cor:best_approximation}
Under Assumption~\ref{assm:basis_normalization}, for every
$V\in L^2(\rho)$,
\[
    \inf_{W\in U_M}\|V-W\|_{\rm E}^2
    =
    \inf_{W\in U_M}\|V-W\|_{L^2(\rho)}^2.
\]
If in addition $DV\in L^2(\rho';\mathbb R^d)$, then
\[
    \inf_{W\in U_M}\|V-W\|_{\rm EF}^2
    =
    \inf_{W\in U_M}
    \Bigl[
        w_{\rm E}\|V-W\|_{L^2(\rho)}^2
        +w_{\rm F}\|D(V-W)\|_{L^2(\rho';\mathbb R^d)}^2
    \Bigr].
\]
\end{corollary}

Thus the best-approximation error for energy observations is the
$L^2(\rho)$ error, and for combined observations a weighted Sobolev-type
error in which the function is measured in $L^2(\rho)$ and its derivative
in $L^2(\rho')$; when $\rho'=\rho$ this is a weighted $H^1(\rho)$ error
(cf.~Remark~\ref{rem:origin_exchange}).

Assumption~\ref{assm:basis_normalization} also ensures that averaging
preserves orthonormality. Define the population energy and force Gram
matrices
\begin{align*}
    \bigl(\mathbb E[G^{\rm E}]\bigr)_{mn}
    &:=
    \mathbb E_{\mathcal Z\sim\rho^{\otimes K}}
    \left[\Phi_m^{\rm E}(\mathcal Z)\Phi_n^{\rm E}(\mathcal Z)\right],
    \\
    \bigl(\mathbb E[G^{\rm F}]\bigr)_{mn}
    &:=
    \mathbb E_{\mathcal Z'\sim\rho'^{\otimes K'}}
    \left[\Phi_m^{\rm F}(\mathcal Z')\cdot\Phi_n^{\rm F}(\mathcal Z')\right],
\end{align*}
and $\mathbb E[G^{\rm EF}]:=w_{\rm E}\mathbb E[G^{\rm E}]+w_{\rm F}\mathbb E[G^{\rm F}]$,
which is the expectation of the empirical Gram matrix
\eqref{eq:empirical_gram}.

\begin{lemma}[Population Gram matrix]
\label{lem:population_gram}
Under Assumptions~\ref{assm:indep} and~\ref{assm:basis_normalization},
the energy features $\{\Phi_m^{\rm E}\}_{m=1}^M$ are orthonormal in
$L^2(\rho^{\otimes K})$, that is, $\mathbb E[G^{\rm E}]=I_M$, and
\begin{equation}\label{eq:population_gram_identities}
    \bigl(\mathbb E[G^{\rm F}]\bigr)_{mn}
    =
    \int_{(\mathbb R^d)^\nu}
        D\varphi_m(\bm z)\cdot D\varphi_n(\bm z)\,d\rho'(\bm z).
\end{equation}
\end{lemma}

\begin{proof}
Write $\mathcal Z=\{\bm z^{(1)},\ldots,\bm z^{(K)}\}$. Separating equal
and distinct cluster indices in
$\Phi_m^{\rm E}(\mathcal Z)\Phi_n^{\rm E}(\mathcal Z)$ and using
independence,
\begin{equation}\label{eq:energy_gram_expansion}
    \mathbb E\!\left[\Phi_m^{\rm E}(\mathcal Z)\Phi_n^{\rm E}(\mathcal Z)\right]
    =
    \int\varphi_m\varphi_n\,d\rho
    +(K-1)\Bigl(\int\varphi_m\,d\rho\Bigr)\Bigl(\int\varphi_n\,d\rho\Bigr).
\end{equation}
By Assumption~\ref{assm:basis_normalization} the right-hand side equals
$\delta_{mn}$: for $m,n\geq2$ by orthonormality and the zero-mean
condition; for $m=1$, $n\geq2$ both terms vanish since
$\int\varphi_1\varphi_n\,d\rho=K^{-1/2}\int\varphi_n\,d\rho=0$; and for
$m=n=1$ it equals $\tfrac1K+\tfrac{K-1}K=1$. Hence $\mathbb E[G^{\rm E}]=I_M$.

For the force part, the same expansion gives
\[
    \mathbb E\!\left[\Phi_m^{\rm F}(\mathcal Z')\cdot\Phi_n^{\rm F}(\mathcal Z')\right]
    =
    \int D\varphi_m\cdot D\varphi_n\,d\rho'
    +(K'-1)\Bigl(\int D\varphi_m\,d\rho'\Bigr)\cdot\Bigl(\int D\varphi_n\,d\rho'\Bigr).
\]
Since $\varphi_m$ is $O(d)$-invariant, $D\varphi_m$ is odd, and inversion
symmetry of $\rho'$ guarantees that $\int D\varphi_m\,d\rho'=0$. The second term
therefore vanishes.
\end{proof}
Following \cite{pmlr-v235-adcock24a}, an observation family $\mathcal O$ is
\emph{nondegenerate} on $U_M$ with respect to a norm $\|\cdot\|_{\mathbb X}$
if there exist constants $0<\alpha\leq\beta<\infty$ such that
\[
    \alpha\|u\|_{\mathbb X}^2
    \leq
    \mathbb E\|\mathcal O(u)\|^2
    \leq
    \beta\|u\|_{\mathbb X}^2
    \qquad\text{for all }u\in U_M,
\]
where for the combined family $\mathbb E\|\mathcal O(u)\|^2$ denotes the
weighted sum of the energy and force terms. The lower bound ensures that no
nonzero direction is invisible to the observations; the upper bound
controls the size of the observation operator.

\begin{corollary}[Observation geometry on $U_M$]
\label{cor:population_gram_geometry}
Under Assumption~\ref{assm:indep}, for every
$u=\sum_{m=1}^M\theta_m\varphi_m\in U_M$,
\begin{align}
    \theta^\top\mathbb E[G^{\rm E}]\theta
    &=
    \mathbb E_{\mathcal Z\sim\rho^{\otimes K}}
    \left|\mathcal O_{\rm E}(u;\mathcal Z)\right|^2
    =
    \|u\|_{\rm E}^2,
    \label{eq:energy_observation_geometry_UM}
    \\
    \theta^\top\mathbb E[G^{\rm EF}]\theta
    &=
    w_{\rm E}\,\mathbb E_{\mathcal Z\sim\rho^{\otimes K}}
    \left|\mathcal O_{\rm E}(u;\mathcal Z)\right|^2
    +w_{\rm F}\,\mathbb E_{\mathcal Z'\sim\rho'^{\otimes K'}}
    \left\|\mathcal O_{\rm F}(u;\mathcal Z')\right\|_2^2
    =
    \|u\|_{\rm EF}^2.
    \label{eq:energy_force_observation_geometry_UM}
\end{align}
Consequently, both observation families are nondegenerate with respect to
their induced norms with $\alpha=\beta=1$; this holds by construction of
the norms.
\end{corollary}
\begin{proof}
By linearity,
$\mathcal O_{\rm E}(u;\mathcal Z)=\sum_m\theta_m\Phi_m^{\rm E}(\mathcal Z)$,
so $\mathbb E|\mathcal O_{\rm E}(u;\mathcal Z)|^2=\theta^\top\mathbb E[G^{\rm E}]\theta$
by definition of the population Gram matrix, and
Lemma~\ref{lem:population_observation_geometry} identifies this with
$\|u\|_{\rm E}^2$. The force term is analogous.
\end{proof} 

The remaining question is how many observations are needed for the
empirical Gram matrix to preserve the population geometry. We address
that question using the variation of the observation operators.

\subsection{Finite-sample recovery}\label{sec:finite_sample_recovery}

We now determine how many observations are sufficient for the empirical
least-squares problem to retain the population geometry identified above.
Following Adcock et al.~\cite{pmlr-v235-adcock24a}, this is governed by the
\emph{variation} of the observation operators on the model space.

Let $\mathcal O_{\mathcal Z}$ denote an observation operator with random
input $\mathcal Z$. The variation of $U_M$ with respect to $\mathcal O_{\mathcal Z}$ is defined as

\begin{equation}\label{eq:general_variation}
    \Lambda(U_M)
    :=
    \operatorname*{ess\,sup}_{\mathcal Z}
    \sup_{\substack{u\in U_M\\\|u\|_{\mathbb X}=1}} \|\mathcal O_{\mathcal Z}(u)\|_2^2.
\end{equation}

Thus, $\Lambda(U_M)$ measures the largest size of a single observation of
a unit-norm model potential. It depends on the space $U_M$, the observation
operator, and the population norm, but not on the choice of basis used to
represent $U_M$.

For energy-only observations, define
\begin{align}
    \Lambda_{\rm E}(U_M)
    =
    \operatorname*{ess\,sup}_{\mathcal Z}
    \sup_{\substack{u\in U_M\\\|u\|_{\rm E}=1}}
    \left|
        \mathcal O_{\rm E}(u;\mathcal Z)
    \right|^2.
\label{eq:energy_variation}
\end{align}

For the combined setting, define 

\begin{align}
    \Lambda_{\rm E}^{\rm EF}(U_M)
    &:=
    \operatorname*{ess\,sup}_{\mathcal Z}
    \sup_{\substack{u\in U_M\\\|u\|_{\rm EF}=1}}
    \left|
        \mathcal O_{\rm E}(u;\mathcal Z)
    \right|^2,
    \label{eq:combined_energy_variation}
    \\
    \Lambda_{\rm F}^{\rm EF}(U_M)
    &:=
    \operatorname*{ess\,sup}_{\mathcal Z'}
    \sup_{\substack{u\in U_M\\\|u\|_{\rm EF}=1}}
    \left\|
        \mathcal O_{\rm F}(u;\mathcal Z')
    \right\|_2^2.
    \label{eq:combined_force_variation}
\end{align}

We now apply the recovery result of
Adcock et al.~\cite{pmlr-v235-adcock24a} to our energy and force
observation families. This gives the following finite-sample guarantee.

\begin{theorem}[Finite-sample recovery]
\label{thm:finite_sample_recovery}

Suppose Assumptions~\ref{assm:fixedK} and
\ref{assm:indep} hold, and let
$\epsilon\in(0,1)$.
\begin{enumerate}[(i)]
    \item Let $V_\theta^{\rm E}\in U_M$ be an exact least-squares
    estimator obtained from $N$ energy observations. There exists a constant $c>0$, such that if
\begin{equation}\label{eq:energy_sample_complexity}
        N
        \geq
        c \,\Lambda_{\rm E}(U_M)
        \log\!\left(\frac{2M}{\epsilon}\right),
\end{equation}
then the truncated estimator $\check V_\theta^{\rm E}
        :=
        \min\left\{
            1,
            \frac{\|V\|_{\rm E}}{\|V_\theta^{\rm E}\|_{\rm E}}
        \right\}
        V_\theta^{\rm E}$
    satisfies
    \begin{equation}\label{eq:energy_recovery_bound}
        \mathbb E
        \left\|V- \check V_\theta^{\rm E}\right\|_{\rm E}^2
        \lesssim
        \inf_{W\in U_M}
        \|V-W\|_{\rm E}^2
        +
        C\epsilon \|V \|_{\rm E}^2.
    \end{equation}

    \item Suppose $w_{\rm E},w_{\rm F}>0$, and let $V_\theta^{\rm EF}\in U_M$ be an exact least-squares
    estimator obtained from $N$ energy observations and $N'$ force
    observations. There exists a constant $c>0$, such that if 
    \begin{equation}
    N
    \geq c w_{\rm E} \Lambda_{\rm E}^{\rm EF}(U_M)
    \log\!\left({\textstyle \frac{2M}{\epsilon}}\right)
    \quad \text{and} \quad
    N' \geq c w_{\rm F} 
    \Lambda_{\rm F}^{\rm EF}(U_M)
    \log\!\left({\textstyle \frac{2M}{\epsilon}}\right),
    \end{equation}
    then the truncated estimator $\check V_\theta^{\rm EF}
        :=
        \min\left\{
            1,
            \frac{\|V\|_{\rm EF}}{\|V\|_{\rm EF}}
        \right\}
        V_\theta^{\rm EF}$
    satisfies
    \begin{align}
        \mathbb E
        \left\|V- \check V_\theta^{\rm EF}\right\|_{\rm EF}^2
        \lesssim
        \inf_{W\in U_M}
        \|V-W\|_{\rm EF}^2
        +
        C\epsilon \|V_\theta^{\rm EF}\|_{\rm EF}^2.
        \label{eq:combined_recovery_bound}
    \end{align}
\end{enumerate}
\end{theorem}

\begin{proof}
    By Corollary~\ref{cor:population_gram_geometry}, the energy and force observations are nondegenerate in their respective population norms
with $\alpha=\beta=1$. Since $U_M$ is an $M$-dimensional subspace, applying the recovery
theorem of Adcock et al.~\cite[Thm. 3.3]{pmlr-v235-adcock24a} therefore gives the
sample conditions and the corresponding recovery estimates. 

For (ii), write the loss as a single average over the $N + N'$ observations with scaled weights $\bar{w}_{\rm E} = (N+N') w_{\rm E}/N$ and $\bar{w}_{\rm F} = (N+N')w_{\rm F}/N'$, so that the sampling operators are $\sqrt{\bar{w}_{\rm E}} \mathcal{O}_{\rm E}$ and $\sqrt{\bar{w}_{\rm F}} \mathcal{O}_{\rm F}$. The multimodal form of \cite[Thm. 3.3]{pmlr-v235-adcock24a} then requires $N + N' \geq c \max\{\bar{w}_{\rm E}\Lambda_{\rm E}^{\rm EF}, \bar{w}_{\rm F} \Lambda_{\rm F}^{\rm EF}\} \log(2M/\epsilon)$, which is equivalent to the two stated conditions.
\end{proof}

Here $\epsilon$ is an upper bound on the probability that empirical
nondegeneracy fails. Truncation controls the error on this event, resulting
in the term $C\epsilon\|V\|_{\mathbb X}^2$ in the recovery bound.

\begin{rem}[Best-approximation]
The bounds in~\cref{thm:finite_sample_recovery} are stated in the observation-induced norms
$\|\cdot\|_{\rm E}$ and $\|\cdot\|_{\rm EF}$, since these are the norms
that the population loss naturally controls. Corollary~\ref{cor:best_approximation}
shows, however, that the best-approximation errors in these norms coincide
with those in standard function-space norms. The same recovery bounds hold, respectively, with the best-approximation terms measured in $L^2(\rho)$ and the weighted Sobolev norm.
\end{rem}

In the remainder of this section we study the variations of the energy and force observation families in terms of the numbers of clusters $K$ and $K'$ involved.

\subsubsection{Worst-case sample-complexity bounds}\label{sec:worst_case_bounds}

We first consider the worst-case bounds obtained from the essential-supremum
variations in Theorem~\ref{thm:finite_sample_recovery}.

\begin{proposition}[Generic variation bounds]
\label{prop:generic_variation}
Under Assumptions~\ref{assm:indep} and
\ref{assm:basis_normalization}, there exist constants
$C_{\rm E}(U_M)$ and $C_{\rm D}(U_M)$, independent of $K$ and $K'$,
such that
\[
    \Lambda_{\rm E}(U_M)
    \leq
    1+K C_{\rm E}(U_M),
\]
and 
\[
    \Lambda_{\rm E}^{\rm EF}(U_M)
    \leq
    \frac{1}{w_{\rm E}}\Lambda_{\rm E}(U_M),
    \qquad
\Lambda_{\rm F}^{\rm EF}(U_M)
    \leq
    \frac{K'}{w_{\rm F}}C_{\rm D}(U_M).
\]
\end{proposition}

The proof is given in~\cref{app:proofs_variation_bounds}, where the constants
$C_{\rm E}$ and $C_{\rm D}$ are explicitly defined in~\cref{eq:ce,eq:cd}.

Combining Proposition~\ref{prop:generic_variation} with
Theorem~\ref{thm:finite_sample_recovery}, the weights cancel and the
sufficient sample sizes become
\[
    N \geq c\,\bigl(1+K C_{\rm E}(U_M)\bigr)\log\!\left(\tfrac{2M}{\epsilon}\right),
    \qquad
    N' \geq c\,K' C_{\rm D}(U_M)\log\!\left(\tfrac{2M}{\epsilon}\right).
\]
The weights $w_{\rm E},w_{\rm F}$ determine the norm in which the error
is measured, but not the number of observations required. For a fixed
approximation space, the worst-case sufficient sample sizes scale linearly
with $K$ and $K'$. 
This linear growth is not an artefact of the proof; one can readily check that 
it is a sharp bound. However, these bounds represent worst-case bounds and overestimate the sample sizes observed in practice as we will see next. 

\subsubsection{Expected sample-complexity bounds}
\label{sec:expected_variation}
Replacing the essential supremum in~\cref{eq:general_variation} by an expectation gives the typical observation size~\cite{cohen2017optimal} and, when the empirical quantities concentrate around it, suggests the sample scale seen in practice.
\begin{proposition}[Bounds on the expected variation]
\label{prop:expected_observation_size}
Under Assumptions~\ref{assm:indep} and
\ref{assm:basis_normalization},
\begin{equation}\label{eq:expected_energy_observation_size}
    \mathbb E_{\mathcal Z\sim\rho_K}
    \left[
        \sup_{\substack{u\in U_M\\\|u\|_{\rm E}=1}}
        \left|
            \mathcal O_{\rm E}(u;\mathcal Z)
        \right|^2
    \right]
    =
    M.
\end{equation}
Moreover,
\begin{equation}\label{eq:expected_force_observation_size}
    \mathbb E_{\mathcal Z'\sim\rho_{K'}}
    \left[
        \sup_{\substack{u\in U_M\\\|u\|_{\rm EF}=1}}
        \left\|
            \mathcal O_{\rm F}(u;\mathcal Z')
        \right\|_2^2
    \right]
    \leq
    \frac{M-1}{w_{\rm F}}.
\end{equation}
\end{proposition}

The proof is given in~\cref{app:expected_observation_size}.
Heuristically, replacing the variations in
Theorem~\ref{thm:finite_sample_recovery} by these expectations gives the
sample conditions
\[
    N \gtrsim M\log\!\left(\tfrac{2M}{\epsilon}\right),
    \qquad
    N' \gtrsim (M-1)\log\!\left(\tfrac{2M}{\epsilon}\right),
\]
independent of $K$ and $K'$; the weights cancel as in
\cref{sec:worst_case_bounds}. Making this heuristic rigorous requires
concentration of the empirical variation around its expectation, which we do not
pursue here; cf. the analogous treatment in~\cite{cohen2017optimal}.

\subsubsection{Specialization to pair interactions}
\label{sec:pair_legendre_variation}

Consider the pairwise interaction case $\nu=1$, where $\bm z\in\mathbb R^d$,
$r=|\bm z|\in[r_{\min},r_{\rm cut}]$, and $U_M$ consists of radial functions
$u({\bm z}) \equiv u(r)$ (allowing some mild abuse of notation). We assume that
$\rho=\rho'$ is isotropic with radial marginal $\rho_r$. Let
$\xi:[r_{\min},r_{\rm cut}]\to[-1,1]$ be a $C^1$ bijection with
$0<c_\xi\leq|\xi'|\leq C_\xi<\infty$, and let $\rho_r$ be the measure whose
pushforward under $\xi$ is the uniform probability measure on $[-1,1]$. With
$P_n$ the Legendre polynomial of degree $n$, the basis
\begin{equation}\label{eq:legendre_basis}
    \varphi_1:=\frac{1}{\sqrt K},
    \qquad
    \varphi_{n+1}(\bm z):=\sqrt{2n+1}\,P_n\bigl(\xi(r)\bigr),
    \qquad
    n=1,\ldots,M-1,
\end{equation}
satisfies Assumption~\ref{assm:basis_normalization}, and $U_M$ is the
space of polynomials of degree at most $M-1$ in $\xi(r)$. This is the
setting of the pairwise inteaction experiments in \cref{sec:numerics}.

Although we use the Legendre basis to derive explicit bounds, the
variation constants depend only on the resulting approximation space
$U_M$, not on the basis used to represent it.

\begin{proposition}[Variation bounds for the Legendre space]
\label{prop:legendre_variation}
For the approximation space generated by
\eqref{eq:legendre_basis}, the energy variation satisfies
\begin{equation}\label{eq:legendre_energy_variation}
    \Lambda_{\rm E}(U_M)
    \leq
    1+K(M^2-1).
\end{equation}
For combined energy and force observations,
\begin{align}
    \Lambda_{\rm E}^{\rm EF}(U_M)
    &\leq
    \frac{1}{w_{\rm E}}
    \left[1+K(M^2-1)\right],
    \label{eq:legendre_combined_energy_variation}
    \\
    \Lambda_{\rm F}^{\rm EF}(U_M)
    &\leq
    \frac{K'}{w_{\rm F}}
    \left(\frac{C_\xi}{c_\xi}\right)^2
    (M-1)^2.
    \label{eq:legendre_force_variation}
\end{align}
\end{proposition}

The proof is given in~\cref{app:pair_wise_variation}.

Inserting these bounds into Theorem~\ref{thm:finite_sample_recovery},
the worst-case sample conditions for the Legendre space read
\[
    N \gtrsim \bigl(1+K(M^2-1)\bigr)\log\!\left(\tfrac{2M}{\epsilon}\right),
    \qquad
    N' \gtrsim K'\Bigl(\frac{C_\xi}{c_\xi}\Bigr)^2(M-1)^2
    \log\!\left(\tfrac{2M}{\epsilon}\right),
\]
whereas the expected variations of
Proposition~\ref{prop:expected_observation_size}, with $\dim D(U_M)=M-1$,
suggest $N\gtrsim M\log(2M/\epsilon)$ and $N'\gtrsim(M-1)\log(2M/\epsilon)$.
The two differ by one factor of $K$ or $K'$ and one factor of $M$. 
The numerical results in Section~\ref{sec:numerics} show that the required sample sizes observed in practice are much closer to the expected-variation estimate than to the worst-case bound.

\subsection{Random numbers of clusters}
\label{sec:random_cluster_counts}
Assumption~\ref{assm:fixedK} was made to simplify the presentation. In practice,
the number of clusters varies between configurations. We now show that our
results extend to this case with minor modifications. Let $\kappa$ and $\kappa'$
denote the random numbers of, respectively, clusters in an energy and force
observation. Assume that $\kappa,\kappa'\geq 1$ almost surely and that,
conditionally on $\kappa=k$ and $\kappa'=k'$, the energy and force cluster
coordinates are i.i.d.\ with laws $\rho$ and $\rho'$. 

The observation operators are normalized using the cluster
counts:
\[
    \mathcal O_{\rm E}(u;\mathcal Z)
    :=
    \frac{1}{\sqrt{\kappa}}
    \sum_{j=1}^{\kappa}u(\bm z^{(j)}),
    \qquad
    \mathcal O_{\rm F}(u;\mathcal Z')
    :=
    \frac{1}{\sqrt{\kappa'}}
    \sum_{j=1}^{\kappa'}Du(\bm z'^{(j)}).
\]

Conditioning on $\kappa=k$ and applying
Lemma~\ref{lem:population_observation_geometry} gives
\[
    \mathbb E\left[
        \left|
            \mathcal O_{\rm E}(u;\mathcal Z)
        \right|^2
        \,\middle|\,
        \kappa=k
    \right]
    =
    \|u\|_{L^2(\rho)}^2
    +(k-1)
    \left|
        \int u\,d\rho
    \right|^2.
\]
The conditional force identity does not depend on $k'$, since inversion
symmetry gives $\int Du\,d\rho'=0$. Averaging over $\kappa$ and $\kappa'$ then gives
\begin{align}
    \mathbb E\left|\mathcal O_{\rm E}(u;\mathcal Z)\right|^2
    &=\|u\|_{L^2(\rho)}^2+(\mathbb E[\kappa]-1)
    \left|\int u\,d\rho\right|^2,
    \label{eq:random_count_energy_geometry}\\
    \mathbb E\left\|\mathcal O_{\rm F}(u;\mathcal Z')\right\|_2^2
    &=\|Du\|_{L^2(\rho';\mathbb R^d)}^2.
    \label{eq:random_count_force_geometry}
\end{align}
Hence the induced norms are
\[
    \|u\|_{\rm E}^2:=\|u\|_{L^2(\rho)}^2
    +(\mathbb E[\kappa]-1)\left|\int u\,d\rho\right|^2,
    \qquad
    \|u\|_{\rm EF}^2:=w_{\rm E}\|u\|_{\rm E}^2
    +w_{\rm F}\|Du\|_{L^2(\rho';\mathbb R^d)}^2.
\]
The constant basis function is now normalized by the expected count,
$\varphi_1\equiv 1/\sqrt{\mathbb E[\kappa]}$, so that
$\Phi_1^{\rm E}(\mathcal Z)=\sqrt{\kappa/\mathbb E[\kappa]}$ and
$\mathbb E|\Phi_1^{\rm E}|^2=1$. If $\varphi_2,\dots,\varphi_M$ are
mean-zero and orthonormal in $L^2(\rho)$, then $\mathbb E[G^{\rm E}]=I_M$
and $(\mathbb E[G^{\rm F}])_{mn}=\int D\varphi_m\cdot D\varphi_n\,d\rho'$.
Thus the population Gram identities, the constants $\alpha=\beta=1$, and the best-approximation results remain unchanged
with respect to the modified induced norms.

To obtain worst-case variation bounds, suppose additionally that
\[
    1\leq\kappa\leq K_{\max},
    \qquad
    1\leq\kappa'\leq K'_{\max}.
\]
Repeating the proof of Proposition~\ref{prop:generic_variation} gives
\begin{align}
    \Lambda_{\rm E}(U_M)
    &\leq
    \frac{K_{\max}}{\bar K}
    +
    K_{\max}C_{\rm E}(U_M),
    \label{eq:random_count_energy_variation}
    \\
    \Lambda_{\rm E}^{\rm EF}(U_M)
    &\leq
    \frac{1}{w_{\rm E}}
    \left[
        \frac{K_{\max}}{\bar K}
        +
        K_{\max}C_{\rm E}(U_M)
    \right],
    \label{eq:random_count_combined_energy_variation}
    \\
    \Lambda_{\rm F}^{\rm EF}(U_M)
    &\leq
    \frac{K'_{\max}}{w_{\rm F}}
    C_{\rm D}(U_M).
    \label{eq:random_count_force_variation}
\end{align}
Thus, the population geometry depends on the expected count $\mathbb E[\kappa]$,
whereas the worst-case variation bounds depend on $K_{\max}$ and
$K'_{\max}$. The expected observation-size results are unchanged.
Theorem~\ref{thm:finite_sample_recovery} therefore holds without
Assumption~\ref{assm:fixedK}, with the induced norms and variation bounds
replaced by those above.

\section{Numerical Results}\label{sec:numerics} 
In this section, we illustrate the theoretical results through numerical
experiments. We examine the geometry induced by energy and force observations,
the conditioning of the resulting least-squares problems, and the dependence
of the recovery error on the number of training samples.
\subsection{Pairwise Setting}
To illustrate the theoretical results of Section~\ref{sec:theory}, we initially consider the problem of learning a pair potential from averaged observations. Throughout, we work with the Lennard-Jones potential,
\begin{equation}
    V(r) = 4\epsilon \left[ \left(\frac{\sigma}{r}\right)^{12} - \left(\frac{\sigma}{r}\right)^6 \right].
\end{equation}
We approximate $V$ using a model class $\mathbb{U}_M = \mathrm{span}\{\varphi_1,\dots,\varphi_M\}$ and study recovery from energy and force observations as described in Section~\ref{sec:theory}. 

In the numerical experiments, the basis functions are constructed after applying a change of variables to the pair distance \(r \in [r_{\rm min}, r_{\rm max}]\): 
\[
\xi(r)
=
A_0 + A_1 \left(
1+
\alpha
\frac{(r/r_0)^q}{1+(r/r_0)^{q-p}}
\right)^{-1},
\]
where $A_0, A_1$ are chosen such that the range of $\xi$ becomes the interval $[-1,1]$. 

The basis functions are then taken to be polynomials in the transformed coordinate $\xi(r)$, 
\[
\varphi_m(r)=P_m(\xi(r)),
\]
for example the Legendre polynomials or orthogonal polynomial adapted to  the data distribution. Unless otherwise stated, we use
\[
\alpha=1,\qquad p=2,\qquad q=4,\qquad r_0=r_m.
\]
This construction is chosen to mimic the radial basis construction used in practice in ACE potentials to concentrate resolution near the
equilibrium bond length $r_0$ while reducing resolution near $r=0$ and
near the cutoff $r_{\rm cut}$~\cite{Drautz2019,DUSSON2022110946,10.1063/5.0158783}. Throughout, we take $w_{\rm E}=w_{\rm F}=1$.

\subsubsection{Matched Sampling and Basis Orthogonality}\label{sec:matched}
In the first setting, the sampling distribution is chosen to match the orthogonality measure of the basis functions. More precisely, let $x(r)$ denote the coordinate transformation used to define the basis, and let $\rho_x$ denote the corresponding orthogonality measure in the transformed coordinate. The pair distances are then sampled according to the distribution $r \sim x^{-1}(\rho_x)$, so that the transformed variable $x(r)$ is distributed according to the orthogonality measure of the basis.

\begin{figure}[htb!]
    \centering
    \begin{subfigure}[b]{0.38\textwidth}
        \centering
\includegraphics[width=\textwidth]{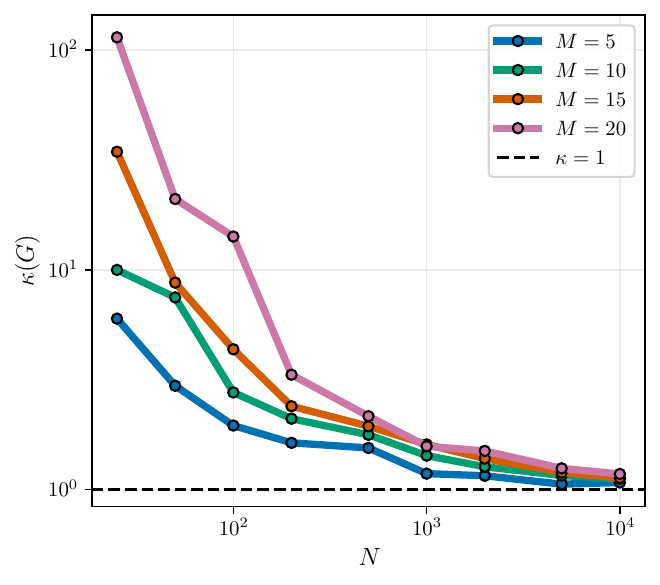}
        \caption{Energy only: $\kappa(G)$ vs.\ $N$ for varying $M$ with $K=18$.}
\label{fig:cond_vs_M_energy}
    \end{subfigure}
\hspace{0.06\textwidth}
    \begin{subfigure}[b]{0.38\textwidth}
        \centering
\includegraphics[width=\textwidth]{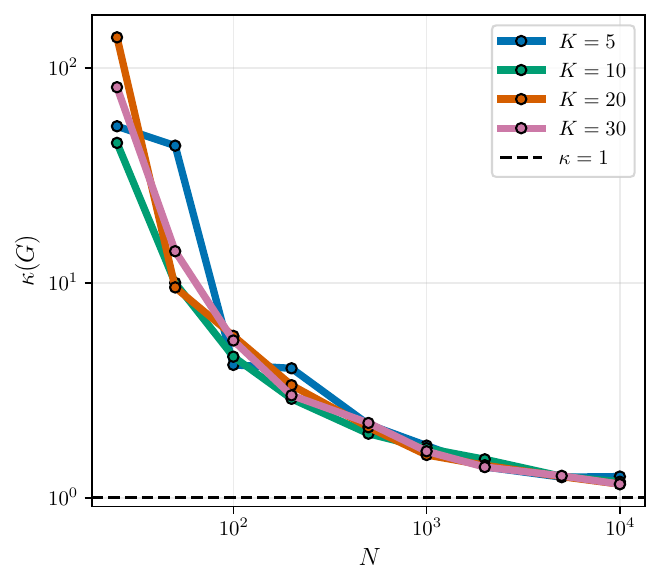}
        \caption{Energy only: $\kappa(G)$ vs.\ $N$ for varying $K$ with $M=18$.}
        \label{fig:cond_vs_K_energy}
    \end{subfigure}
    \begin{subfigure}[b]{0.38\textwidth}
        \centering
\includegraphics[width=\textwidth]{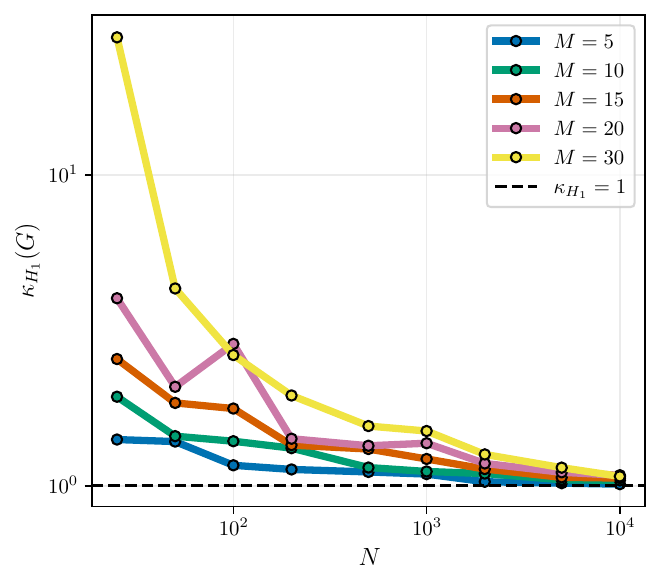}
        \caption{Combined energy and forces: $\kappa(G)$ vs.\ $N$ for varying $M$ with $K=18$.}
\label{fig:cond_vs_M_force}
    \end{subfigure}
\hspace{0.06\textwidth}
    \begin{subfigure}[b]{0.38\textwidth}
        \centering
    \includegraphics[width=\textwidth]{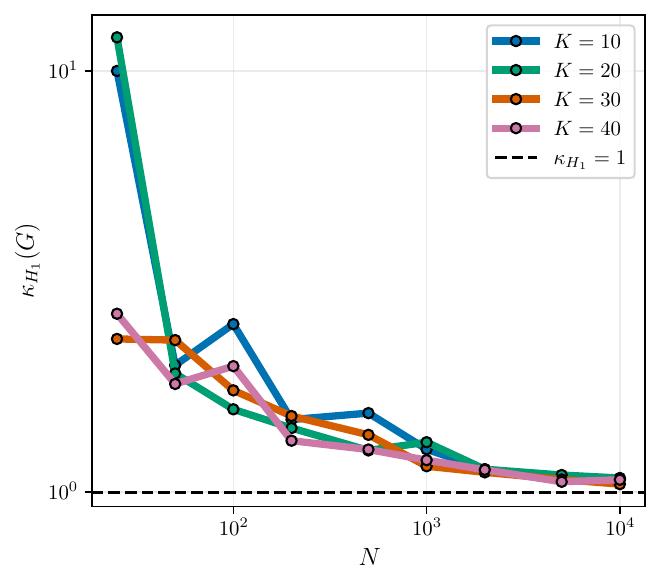}
        \caption{Combined energy and forces: $\kappa(G)$ vs.\ $N$ for varying $K$ with $M=18$.}
        \label{fig:cond_vs_K_force}
    \end{subfigure}

    \caption{Conditioning of the empirical Gram matrix $G$. \textit{Top row:} energy-only setting. \textit{Bottom row:} combined energy and forces setting. In both cases, the condition number decreases and stabilizes close to 1 as $N$ increases, indicating concentration toward a well-conditioned population Gram matrix.}
    \label{fig:gram_conditioning_sec1}
\end{figure}

\paragraph{Gram Matrix Conditioning}

We first examine the conditioning of the empirical Gram matrix for energy-only observations in the geometry induced by the energy observations. As discussed in Section~\ref{sec:theory}, its conditioning determines the stability of the recovery problem.
Given a training set $\{R_i\}_{i=1}^N$, we consider the empirical Gram matrix $
G_N = \frac{1}{N} \sum_{i=1}^N \varphi(R_i)\varphi(R_i)^\top,$
where $\varphi(R)$ denotes the vector of averaged basis functions. 

Figures~\ref{fig:cond_vs_M_energy} and~\ref{fig:cond_vs_K_energy} show the condition numbers $\kappa(G_N)$ as a function of the number of samples $N$, for varying values of the basis size $M$ (left) and the number of pairs per configuration $K$ (right). For small $N$, we observe that the Gram matrix is poorly conditioned, particularly for larger values of $M$. As $N$ increases, the condition number decreases and stabilizes, suggesting concentration of the empirical Gram matrix around its population. The condition numbers for $K=10$ and $K=30$ are comparable, consistent with the theory.

\paragraph{$H^1_D$-scaled condition number}

To compare the conditioning of the empirical regression problem in a way that is consistent with the target $H^1_D$ error, we do not compute the condition number of the empirical Gram matrix directly. Instead, we compute a condition number relative to the $H^1_D$ metric. This is equivalent to using an $H^1_D$ orthogonal basis.

Let $G_{H^1}$ be defined as 
\[
    (G_{H^1})_{mn}
    :=
    \int \varphi_m(r)\varphi_n(r)\,d\rho(r)
    +
    \int \varphi_m'(r)\varphi_n'(r)\,d\rho(r).
\]
then $\kappa_{H^1}(G) :=
    \kappa(
        G_{H^1}^{-1/2}G G_{H^1}^{-1/2}
    )$, where $\kappa$ denotes the standard $\ell^2$-condition number.

This measures how well the empirical Gram matrix $G$
approximates the target $H^1$ metric $G_{H^1}$. If $G \approx G_{H^1}$, then
$G_{H^1}^{-1/2}GG_{H^1}^{-1/2}\approx I$,
and consequently
$\kappa_{H^1}(G)\approx 1$.
In this regime, $G$ defines a norm on the approximation space that is equivalent to the $H^1$ norm. If $G$ is positive definite and $\kappa_{H^1}(G)$ remains bounded as $M$ increases, no direction in the approximation space becomes poorly observed in the $H^1$ norm.
Figures~\ref{fig:cond_vs_M_force} and~\ref{fig:cond_vs_K_force} show
$\kappa_{H^1}(G_N)$ as a function of the number of samples $N$, for varying
values of the basis size $M$ (left) and the number of pairs per
configuration $K$ (right). In all cases, $\kappa_{H^1}(G_N)$ decreases
toward $1$ as $N$ increases, indicating convergence of the empirical Gram
matrix to the target $H^1_D$ metric which indicates $H^1_D$ identifiability is recovered as the sample size grows.

\begin{figure}[t]
    \centering
\includegraphics[width=0.8\textwidth]{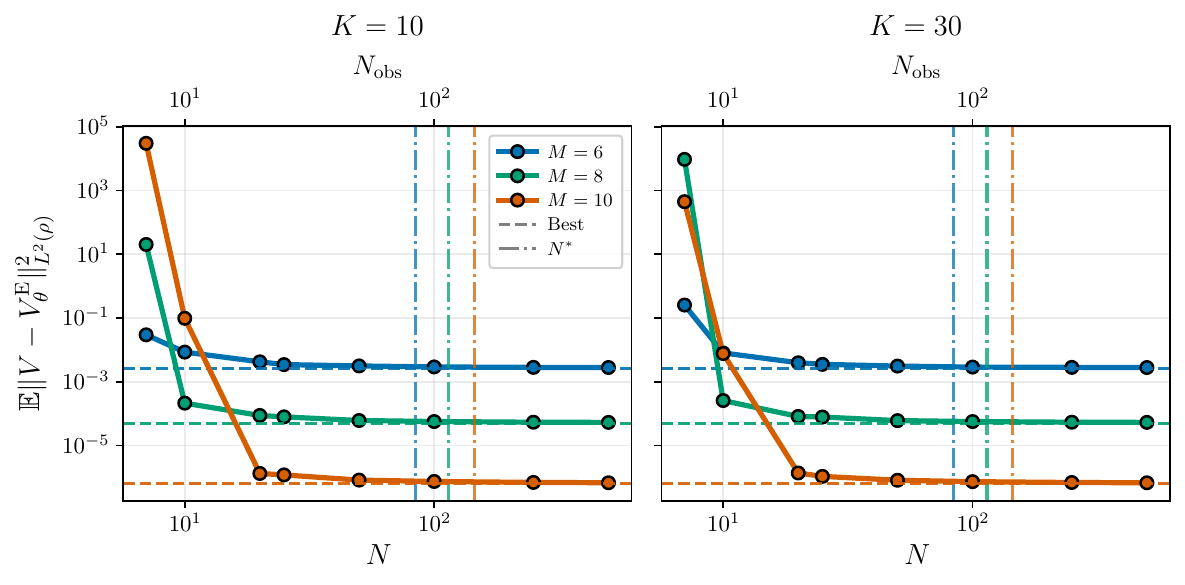}
    \caption{Reconstruction error $\|V - \hat V_\theta\|_{L^2(\rho)}^2$ as a function of the number of samples $N$ for energy-only observations. Dashed horizontal lines indicate the approximation error and the vertical lines indicate the theoretical sample complexity needed to achieve the recovery guarantees.}
    \label{fig:l2_energy_noise}
\end{figure}
\paragraph{Recovery in $L^2$ and $H^1$} We now examine recovery of the pair potential from energy only observations. In this setting, we measure reconstruction error in the $L^2(\rho)$ norm.
Figure~\ref{fig:l2_energy_noise} shows the empirical error $\|V - \hat V_\theta\|_{L^2(\rho)}^2$ as a function of the number of samples $N$, for two values of pairs per configuration, $K=10$ and $K=30$, and varying basis sizes $M$.
In all cases, the error decreases rapidly for small $N$ before approaching the best-approximation error. Increasing $M$ lowers this error but requires more samples to reach it.

\begin{figure}[t]
\centering
\includegraphics[width=0.8\textwidth]{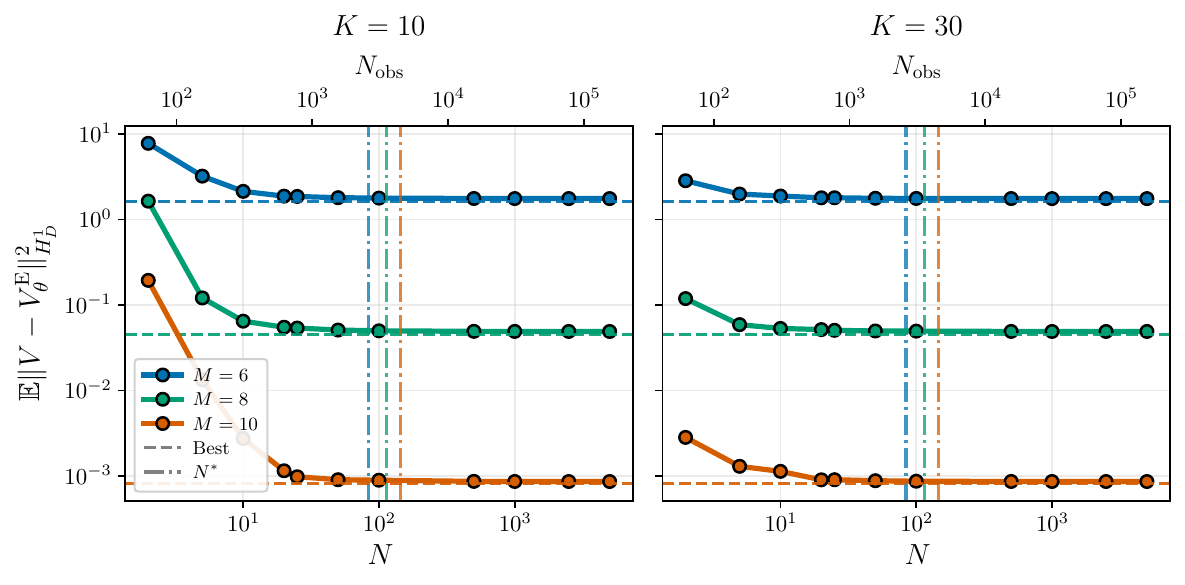}
     \caption{Reconstruction error $\|V - \hat V_\theta\|_{H^1_D}^2$ as a function of the number of samples $N$ for combined energy and force observations.}
    \label{fig:H1_EF}
\end{figure}

To study the recovery from combined energy and force observations, we measure the reconstruction error in the $H^1_D$ norm,
which incorporates both the function and its derivative.  
Figure~\ref{fig:H1_EF} shows the empirical error $\|V - \hat V_\theta\|_{H^1_D}^2$ as a function of the number of samples $N$, for $K=10$ and $K=30$ and varying basis sizes $M$.
As in the energy-only setting, the error decreases rapidly with increasing
sample size before approaching the approximation-error plateau. 
The results for $K=10$ and $K=30$ are analogous. This is consistent with the expected-variation bounds which are independent of $K$.

\FloatBarrier

\subsubsection{Mismatched Sampling}\label{sec:mismatched}

We now verify that changing the basis used to represent a fixed
approximation space does not affect the conditioning of the least-squares solution (assuming exact arithmetic). 
We draw the pair distances from the distribution
\[
    \rho(r)\propto r^2,
    \qquad
    r\in[r_{\min},r_{\rm cut}],
\]
and compare two bases of the same space $U_M$. The first is the Legendre
basis defined in~\eqref{eq:legendre_basis}. Since this basis is not
orthonormal with respect to the observation geometry induced by $\rho$,
we construct a second basis by applying a coordinate transformation to the original basis.

For energy observations, define $
    \Phi^{\rm E}(\mathcal Z)
    :=
    \bigl(
        \Phi_1^{\rm E}(\mathcal Z),\ldots,
        \Phi_M^{\rm E}(\mathcal Z)
    \bigr),$ and let $G_{\rm E}$ be the population Gram matrix of the original basis. Define the transformed energy features by
\[
    \widetilde{\Phi}^{\rm E}(\mathcal Z)
    :=
    \Phi^{\rm E}(\mathcal Z)G_{\rm E}^{-1/2},
\]
so that
\[
    \mathbb E\left[
        (\widetilde{\Phi}^{\rm E})^\top
        \widetilde{\Phi}^{\rm E}
    \right]
    =I_M.
\]
We apply the same transformation using the combined energy-force Gram matrix $G_{\rm EF}$. Since both $G_{\rm E}^{-1/2}$ and $G_{\rm EF}^{-1/2}$ are invertible, these transformations only change the coordinates used to represent $U_M$, the approximation space itself remains unchanged.

\begin{figure}[!htb]
    \centering
    \begin{subfigure}[b]{\textwidth}
        \centering
        \includegraphics[width=0.9\textwidth]{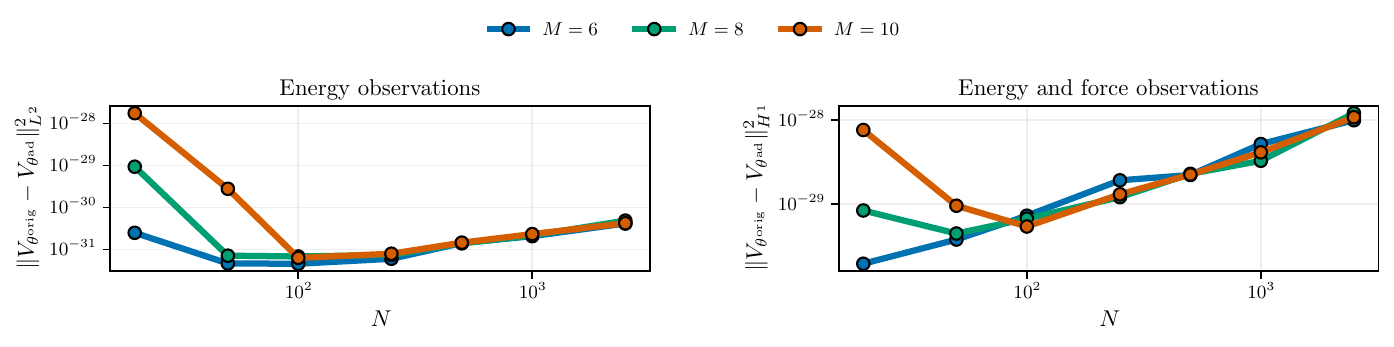}
        \caption{$K=10$.}
    \end{subfigure}
    \vspace{1em}
    \begin{subfigure}[b]{\textwidth}
        \centering
        \includegraphics[width=0.9\textwidth]{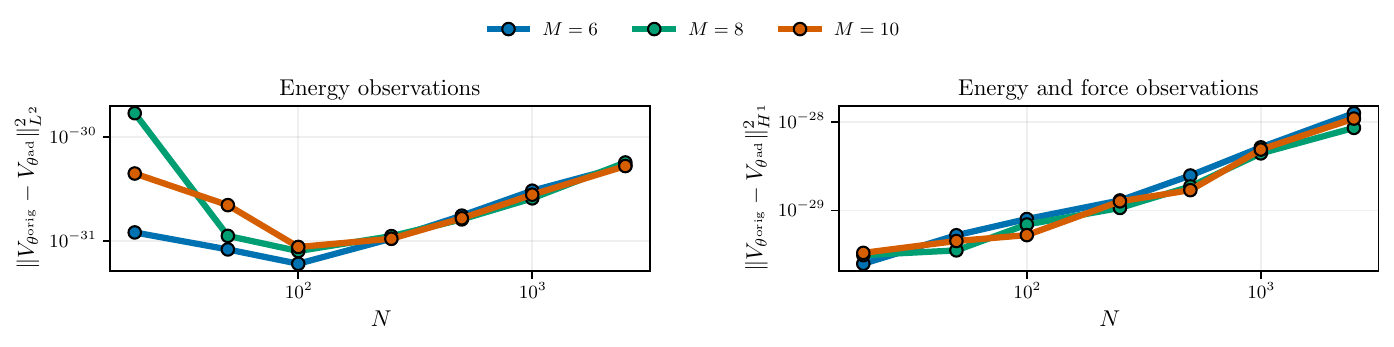}
        \caption{$K=30$.}
    \end{subfigure}
    \caption{Difference between the potentials fitted in the original and
population-adapted bases for $K=10$ (a) and $K=30$ (b), confirming that
the two bases produce the same fitted potential.}
\label{fig:fitted_function_difference}
\end{figure}

\begin{figure}[htb]
    \centering
    \begin{subfigure}[b]{0.7\textwidth}
        \centering
        \includegraphics[width=\textwidth]{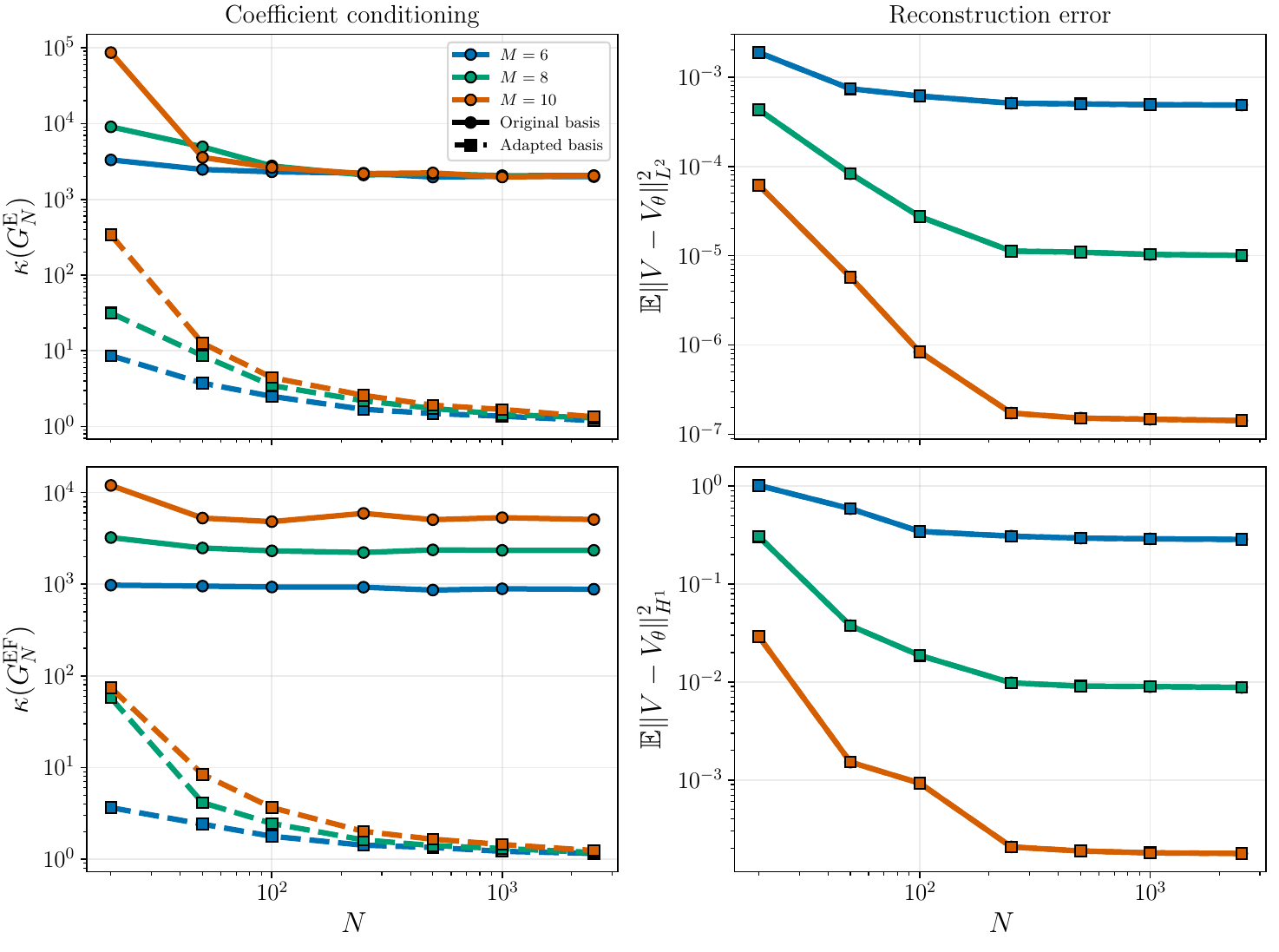}
        \caption{$K=10$.}
        \label{fig:basis_invariance_K10}
    \end{subfigure}
    \vspace{1em}
    \begin{subfigure}[b]{0.7\textwidth}
        \centering
        \includegraphics[width=\textwidth]{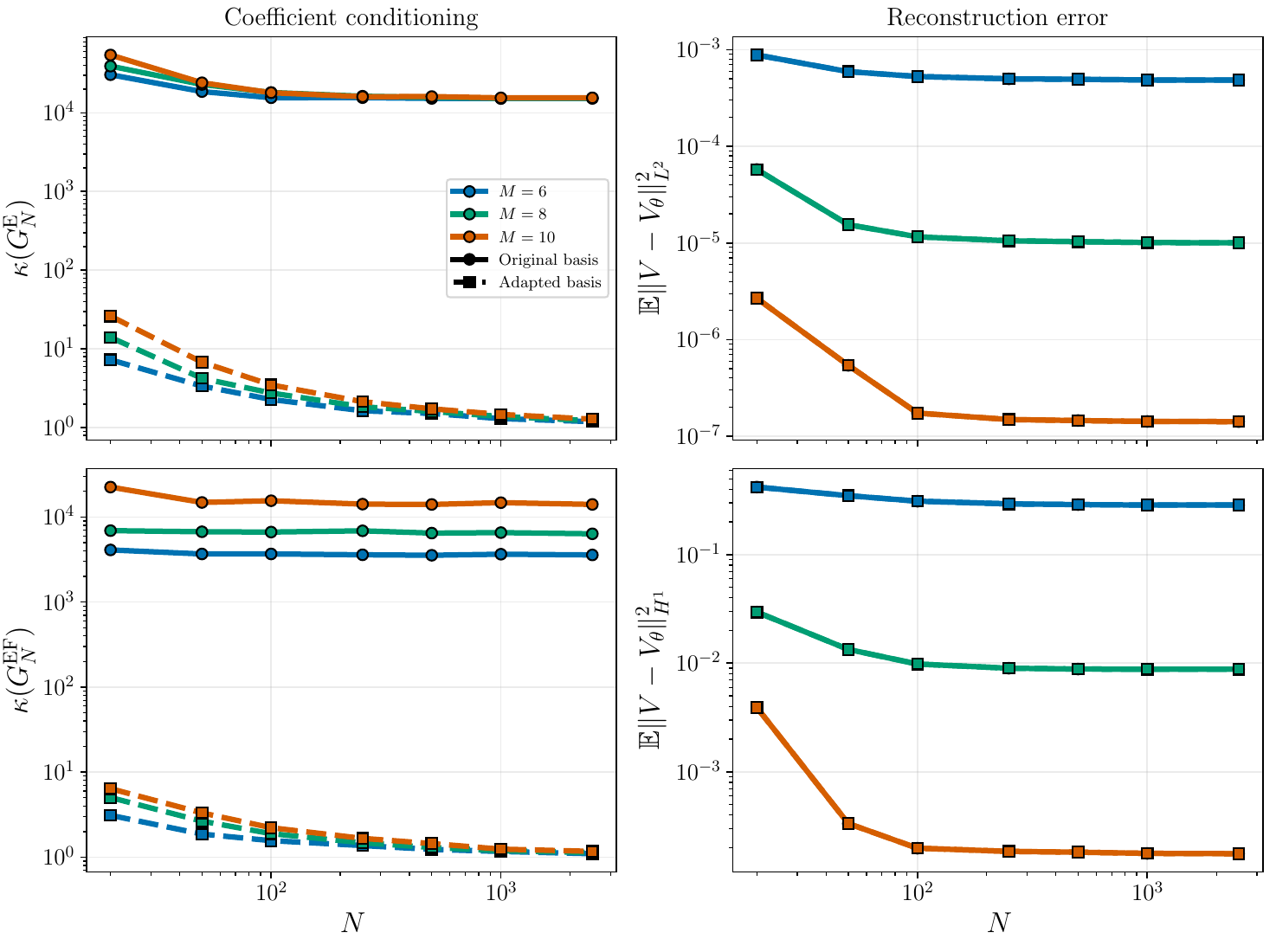}
        \caption{$K=30$.}
        \label{fig:basis_invariance_K30}
    \end{subfigure}
    \caption{Conditioning and reconstruction errors for the original and
population-adapted bases with $K=10$ (a) and $K=30$ (b). The top and
bottom rows show energy-only and combined energy-force recovery,
respectively. The adapted basis improves the conditioning of
the empirical Gram matrix, while the reconstruction error curves
for the two bases lie on top of each other, since they span the same
approximation space.}
    \label{fig:basis_invariance_sec2}
\end{figure}

Figures~\ref{fig:basis_invariance_K10}
and~\ref{fig:basis_invariance_K30} compare the two representations for
$K=10$ and $K=30$, respectively. The population-adapted coordinates
reduce the condition numbers of the empirical Gram matrices by several
orders of magnitude. Despite this difference in conditioning, the reconstruction
errors obtained from the two bases are indistinguishable. Let
$\theta^{\rm orig}$ and $\theta^{\rm ad}$ denote the coefficients
computed in the original and adapted coordinates, respectively. Then $V_{\theta^{\rm orig}}
    =
    V_{\theta^{\rm ad}}$
up to numerical precision, as shown in
Figure~\ref{fig:fitted_function_difference}. This confirms that the exact least-squares estimator depends on the
approximation space and the observations, but not on the basis used to
represent that space. A better-conditioned basis can improve
the numerical coefficient problem without changing the fitted potential.

\FloatBarrier

\subsubsection{Recovery with a gap in the sampling distribution}
We now keep the approximation
space fixed and instead remove an interval from the support of the training
distribution. This is a situation that routinely occurs when fitting a potential to a crystalline solid dataset. 
Specifically, we take $K=10$, $M=10$, and sample from
\[
    \rho_{\rm gap}(r)
    \propto
    r^2
    \mathbf 1_{
        [r_{\min},\,1.15]\cup[1.65,\,r_{\rm cut}]
    }(r).
\]
Thus, the training data contain no pair distances in $I_{\rm gap}:=(1.15,1.65).$

\begin{figure}[!htb]
    \centering
    \includegraphics[width=0.8\textwidth]{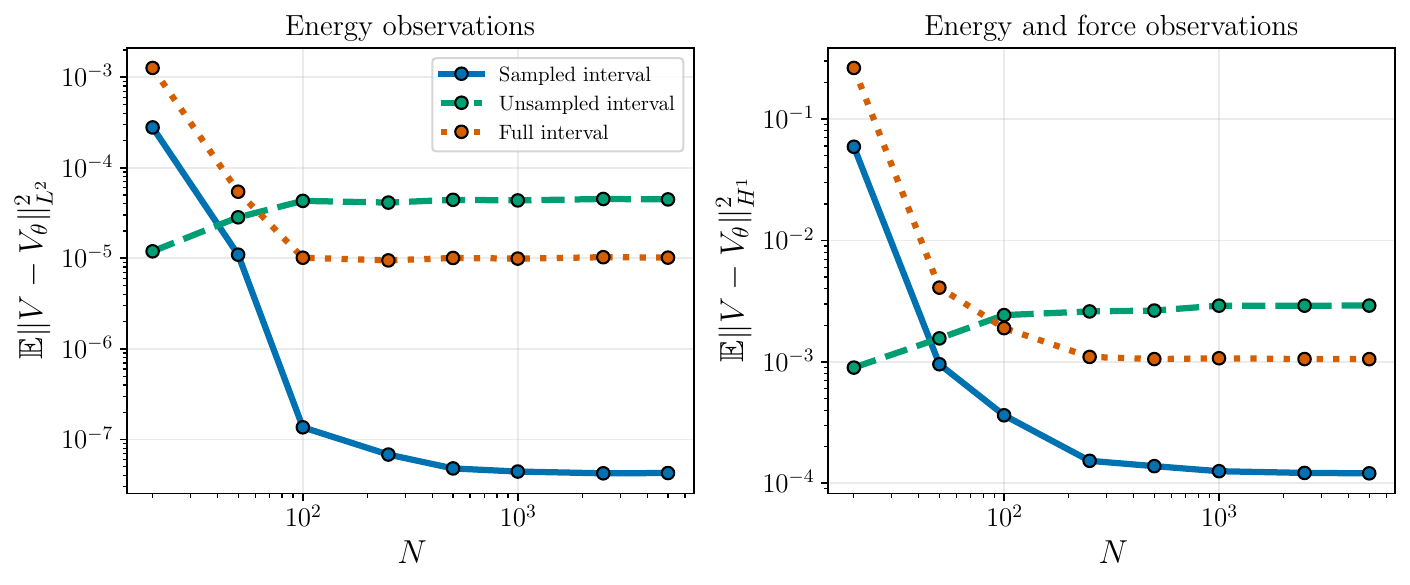}
    \caption{
    Reconstruction errors for a training distribution with no samples in
    $I_{\rm gap}=(1.15,1.65)$. Errors are evaluated on the support of the
    training distribution, inside the unsampled interval, and over the full interval. Increasing the number of samples improves the
    reconstruction on the sampled region but does not reduce the error in
    the gap.
    }
    \label{fig:gap_errors}
\end{figure}

Figure~\ref{fig:gap_errors} shows the reconstruction error separately on
the support of the training distribution, inside the unsampled interval,
and over the full interval $[r_{\min},r_{\rm cut}]$. As the number of
samples increases, the error on the sampled region decreases and approaches
its approximation floor. The error inside $I_{\rm gap}$, however, does
not decrease. The same behaviour occurs for combined energy and force observations.
Derivative information does not provide direct information
inside an interval that is never sampled. Therefore, the population loss controls the potential only on the support of
the sampling distribution and recovery outside this support requires
additional assumptions, prior information, or additional data.

\subsection{Higher Body-Order Experiments}

We now extend the numerical investigation to three-body ($\nu=2$) and
four-body ($\nu=3$) interactions. The goal is to verify that the results shown in the pairwise setting, holds at higher correlation orders.

\subsection{Three-body}

For the three-body experiment, we use the reference potential
\begin{equation}\label{eq:three_body_reference}
    V_{3\rm b}(\bm z_1,\bm z_2)
    =
    e^{-(r_1+r_2)}
    \left(c_{12}+\frac13\right)^2
    \left[
        1+0.2\sin(4r_1)\sin(4r_2)
    \right],
\end{equation}
where
\[
    r_i:=\|\bm z_i\|_2,
    \qquad
    c_{ij}:=
    \frac{\bm z_i\cdot\bm z_j}{r_ir_j}.
\]

The approximation spaces are constructed from the invariant ACE basis
described in Appendix~\ref{sec:ace_appndx}.

Figure~\ref{fig:3body} shows the reconstruction errors for energy-only
and combined energy-force observations. In both cases, the error decreases
with the number of training configurations and approaches the corresponding
best-approximation error. Increasing $M$ lowers this error but
requires more data before the approximation regime is reached. The results
for $K=10$ and $K=30$ are similar, in agreement with the weak dependence
on $K$ suggested by the expected-variation bounds.

\begin{figure}[!htbp]
    \centering
    \begin{subfigure}[b]{0.8\textwidth}
        \centering
\includegraphics[width=\textwidth]{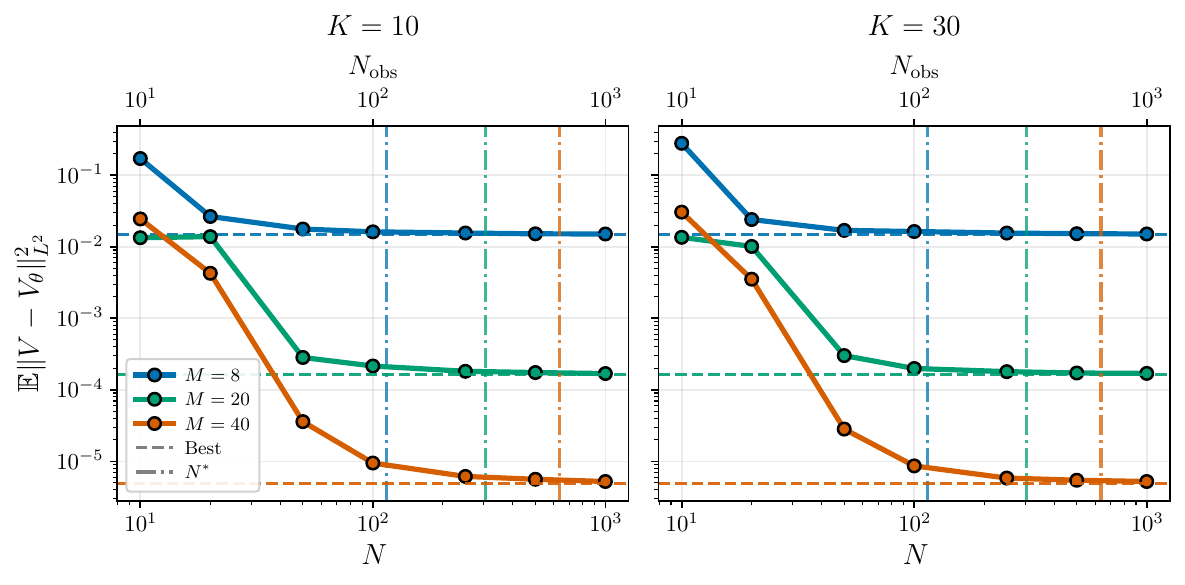}
        \caption{Energy observations}
    \end{subfigure}
    \begin{subfigure}[b]{0.8\textwidth}
        \centering
\includegraphics[width=\textwidth]{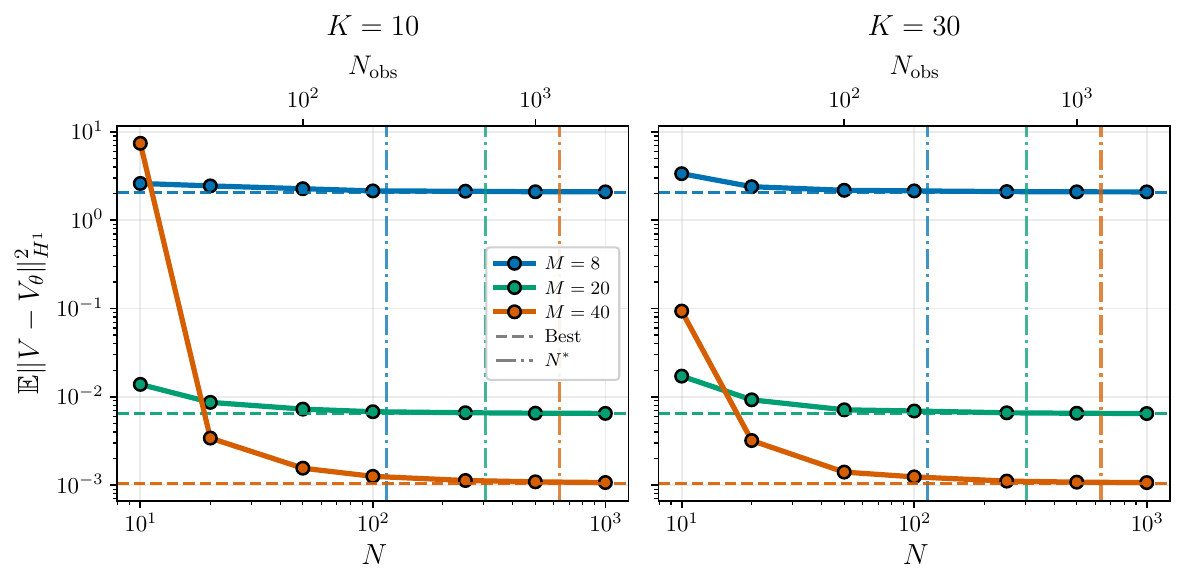}
        \caption{Energy and force observations}
    \end{subfigure}
    \caption{Three-body reconstruction errors for energy-only observations (top row) and
combined energy-force observations (bottom row), with $K=10$ and
$K=30$. Dashed horizontal lines
    indicate the approximation error; vertical lines indicate the
    theoretical sample complexity $N^*$.}
    \label{fig:3body}
\end{figure}

\subsection{Four-body}

For the four-body experiment, we take
\begin{equation}\label{eq:four_body_reference}
\begin{aligned}
V_{4\rm b}(\bm z_1,\bm z_2,\bm z_3)
&=
\exp\!\bigl(-(r_1+r_2+r_3)\bigr)
\left[
    1+\frac12
    \sum_{1\leq i<j\leq3}
    \left(c_{ij}+\frac13\right)^2
\right]
\\
&\quad\times
\left[
    1+
    \frac{0.4}{6}
    \sum_{\substack{i,j,k\in\{1,2,3\}\\
                    i,j,k\ {\rm distinct}}}
    \sin(4r_i)\sin(5r_j)\sin(6r_k)
\right].
\end{aligned}
\end{equation}

\begin{figure}[!htbp]
    \centering
    \begin{subfigure}[b]{0.8\textwidth}
        \centering
\includegraphics[width=\textwidth]{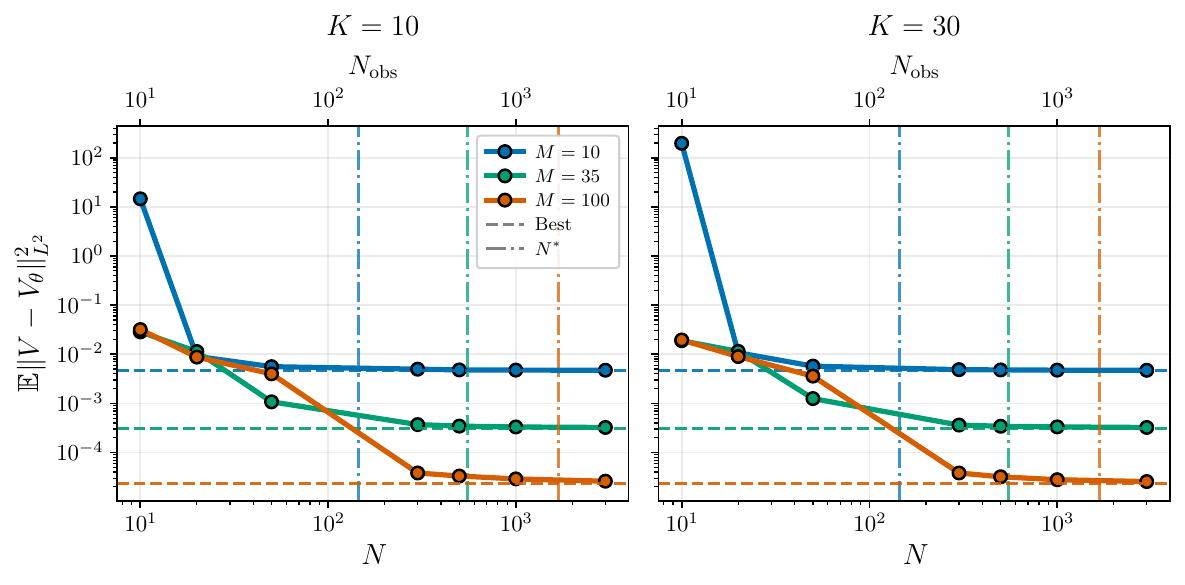}
        \caption{Energy observations}
    \end{subfigure}
    \begin{subfigure}[b]{0.8\textwidth}
        \centering
\includegraphics[width=\textwidth]{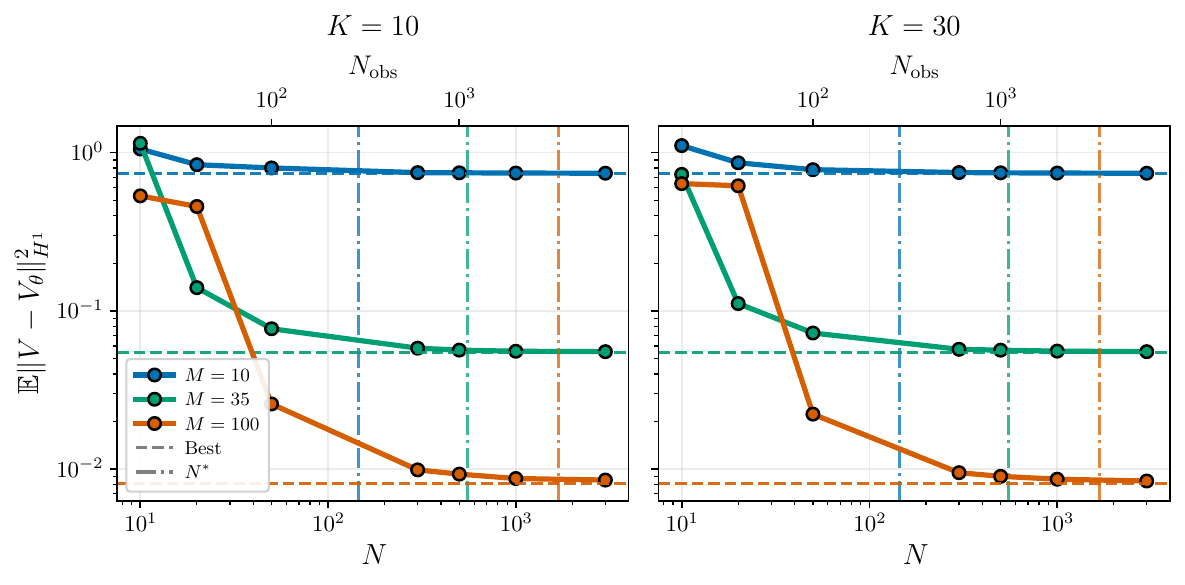}
        \caption{Energy and force observations}
    \end{subfigure}
    \caption{Four-body reconstruction errors for energy-only observations (top row) and
combined energy-force observations (bottom row), with $K=10$ and
$K=30$.}
    \label{fig:4body}
\end{figure}

Figure~\ref{fig:4body} shows the corresponding energy-only and
energy-force reconstruction errors. The same pattern is observed as in
the three-body experiment. The errors approach the best-approximation levels, while increasing $M$ lowers the limiting error but requires more samples.

The three- and four-body experiments therefore have the same basic
recovery behavior as the pairwise experiment. In each case, the reconstruction error approaches
the best error available in the chosen approximation space, and the
observed data requirements are closer to the expected-variation scale shown in~\cref{sec:expected_variation} than
to the worst-case bounds.

\FloatBarrier
\section{Conclusion and Future Work}\label{sec:conclusion}
We introduced an approach to studying machine learning interatomic potentials from the perspective of inverse problems.
The central message of this work is that local interaction models can be recovered from aggregate observations, and that the notion of stable recovery is determined by the type of observations available. 
Energy-only observations identify the potential in an $L^2$ sense, while the inclusion of force information naturally lifts this to $H^1$. 
The data distribution determines weights for the $L^2$ and $H^1$ norms. 
Our approach, employing the framework of \cite{pmlr-v235-adcock24a}, also immediately gives sample complexity bounds. 
Numerical experiments in pairwise, three-body, and four-body settings confirmed the predicted Gram-matrix concentration and validated the observation-induced $H^1$ geometry. A central observation is that the averaging in aggregate observations
does not deteriorate the recovery. In the population limit it only
reweights the constant mode, and while the worst-case sample complexity
bounds grow with the number of clusters per observation, the expected
observation size does not. The numerical tests confirm the latter: only
the number of observations matters, not the number of clusters they
average over.

Our work should be understood as a first step in a highly idealised and simplified setting involving a single correlation order and training data consistent of i.i.d.\ clusters, which enabled an elementary analysis. This establishes a best-case scenario against which results under more realistic assumptions can be measured. Most prominently, it will be necessary to extend our framework to (i) estimating all body orders at the same time, which creates a source of ill-conditioning~\cite{DUSSON2022110946}; and (ii) realistic datasets generated from molecular dynamics type sampling methods where clusters within a structure are highly correlated.

\section*{Acknowledgments} 
The authors acknowledge helpful discussions with Ben Adcock. TT was supported by the Killam Doctoral Scholarship. 

\appendix

\section{Energy and Force Observations}~\label{sec:force_obsv}
Restricting to a correlation order $\nu,$ the total energy is given by \begin{equation}
E(\mathbf{R}) =\sum_{c \in \widetilde{\mathcal{P}}(\mathbf{R})} \widetilde V(c).
\end{equation}

We denote the corresponding collection of relative cluster coordinates by the multiset
\[
    \mathcal Z(\mathbf R)
    :=
    \left\{
        \bm z^{(c)}
        :
        c\in\widetilde{\mathcal P}(\mathbf R)
    \right\}.
\]

Hence, the total energy can be expressed as
\begin{equation}
    E(\mathbf{R}) =
    \sum_{\bm z\in\mathcal{Z}}
    V(\bm z).
\end{equation}

The corresponding atomic forces are obtained by differentiating the total energy with respect to the atomic positions. In particular, the force on atom $i$ is given by:
\begin{align*}
    \bm F_i(\mathbf R)
    &:=
    -\nabla_{\bm r_i}E(\mathbf R)
    \\
    &=
    -\nabla_{\bm r_i}
\left(
\sum_{\substack{c \in \widetilde{\mathcal{P}}(\mathbf{R})\\ \bm{r}_i\in c}}
\widetilde V
\bigl(\bm{r}_i, \bm{r}_{j'_1}, \cdots, \bm{r}_{j'_\nu}
\bigr)
\right)  \\[3pt]
&=
-\nabla_{\bm{r}_i}
\left(
\sum_{\substack{c \in \widetilde{\mathcal{P}}(\mathbf{R})\\ \bm{r}_i\in c}} V
\bigl(
\bm r_{i}-\bm{r}_{j'_1},\,
\bm r_{i}-\bm{r}_{j'_2},\, 
\ldots\, \bm r_{i}-\bm{r}_{j'_\nu}
\bigr)
\right)
 \\[3pt]
&=
\sum_{\substack{c \in \widetilde{\mathcal{P}}(\mathbf{R}) \\ \bm{r}_i\in c}}  \sum_{\ell=1}^{\nu} \nabla_{\bm{r}_{j'_{\ell}}}V
\bigl(
\bm r_{i}-\bm{r}_{j'_1},\,
\bm r_{i}-\bm{r}_{j'_2},\, 
\ldots\, \bm r_{i}-\bm{r}_{j'_\nu}
\bigr)
\end{align*}

Recall the collection of
relative coordinates of the interacting clusters containing $i$ by
\[
    \mathcal Z'(\mathbf R,i)
    :=\left\{
        (\bm r_{j_1}-\bm r_i,\ldots,
         \bm r_{j_\nu}-\bm r_i)
        :
        \begin{array}{l}
            j_1<\cdots<j_\nu,\\
            j_\ell\neq i \quad\text{for all }\ell,\\
            r_{ab}\leq r_{\rm cut}
            \quad\text{for all }
            a,b\in\{i,j_1,\ldots,j_\nu\}
        \end{array}\right\}.
\]
Thus, the force on atom $i$ can be represented as

\begin{equation}
    \bm F_i(\mathbf R)
    :=
    \sum_{\bm{z} \in \mathcal Z'(\mathbf R,i)}
    \sum_{\ell=1}^{\nu}
    \nabla_{\bm z_\ell}V(\bm z)
    \in\mathbb R^d,
\end{equation}
where $\bm z=(\bm z_1,\ldots,\bm z_\nu)\in(\mathbb R^d)^\nu.$

\section{Proofs of the variation bounds}

\subsection{Worst-case (\cref{prop:generic_variation})}\label{app:proofs_variation_bounds}
\begin{proof}
For the energy bound, introduce the pointwise Christoffel function
\begin{equation}
    \mathcal K_{{\rm E},M}(\bm z)
    :=
    \sup_{\substack{
        u\in U_M,\ \int u\,d\rho=0\\
        \|u\|_{L^2(\rho)}=1
    }}
    |u(\bm z)|^2,
\end{equation}
and set
\begin{equation}\label{eq:ce}    C_{\rm E}(U_M)
    :=
    \operatorname*{ess\,sup}_{\bm z\sim\rho}
    \mathcal K_{{\rm E},M}(\bm z).
\end{equation}
Under Assumption~\ref{assm:basis_normalization},
$\varphi_2,\ldots,\varphi_M$ form an orthonormal basis for the mean-zero
functions in $U_M$. Consequently,
\begin{equation}\label{eq:energy_christoffel_features}
    \mathcal K_{{\rm E},M}(\bm z)
    =
    \sum_{m=2}^M|\varphi_m(\bm z)|^2.
\end{equation}
Moreover, since the averaged energy features are orthonormal,
\[
    \Lambda_{\rm E}(U_M)
    =
    \operatorname*{ess\,sup}_{\mathcal Z\sim\rho_K}
    \left[
        1+
        \sum_{m=2}^M
        |\Phi_m^{\rm E}(\mathcal Z)|^2
    \right].
\]
Writing
$\mathcal Z=\{\bm z^{(1)},\ldots,\bm z^{(K)}\}$ and applying
Cauchy-Schwarz gives
\begin{align*}
    \sum_{m=2}^M|\Phi_m^{\rm E}(\mathcal Z)|^2
    &=
    \sum_{m=2}^M
    \left|
        \frac{1}{\sqrt K}
        \sum_{k=1}^K\varphi_m(\bm z^{(k)})
    \right|^2
    \\
    &\leq
    \sum_{k=1}^K
    \sum_{m=2}^M
    |\varphi_m(\bm z^{(k)})|^2
    \\
    &=
    \sum_{k=1}^K
    \mathcal K_{{\rm E},M}(\bm z^{(k)})
    \\
    &\leq
    K C_{\rm E}(U_M).
\end{align*}
Therefore,
\[
    \Lambda_{\rm E}(U_M)
    \leq
    1+K C_{\rm E}(U_M).
\]

The inequality
\[
    \|u\|_{\rm EF}^2
    \geq
    w_{\rm E}\|u\|_{\rm E}^2
\]
also gives
\[
    \frac{
        w_{\rm E}|\mathcal O_{\rm E}(u;\mathcal Z)|^2
    }{
        \|u\|_{\rm EF}^2
    }
    \leq
    \frac{
        |\mathcal O_{\rm E}(u;\mathcal Z)|^2
    }{
        \|u\|_{\rm E}^2
    },
\]
and hence
\[
\Lambda_{\rm E}^{\rm EF}(U_M)
    \leq
    \frac{1}{w_{\rm E}}\Lambda_{\rm E}(U_M).
\]

For the force bound, we introduce
\[
    \mathcal K_{{\rm D},M}(\bm z)
    :=
    \sup_{\substack{
        u\in U_M\\
        \|Du\|_{L^2(\rho';\mathbb R^d)}=1
    }}
    \|Du(\bm z)\|_2^2
\]
and set
\begin{equation}\label{eq:cd}
    C_{\rm D}(U_M)
    :=
    \operatorname*{ess\,sup}_{\bm z\sim\rho'}
    \mathcal K_{{\rm D},M}(\bm z).
\end{equation}
For
$\mathcal Z'=\{\bm z^{(1)},\ldots,\bm z^{(K')}\}$,
Cauchy-Schwarz gives
\begin{align*}
    \left\|
        \mathcal O_{\rm F}(u;\mathcal Z')
    \right\|_2^2
    &=
    \frac{1}{K'}
    \left\|
        \sum_{k=1}^{K'}Du(\bm z^{(k)})
    \right\|_2^2
    \\
    &\leq
    \sum_{k=1}^{K'}\|Du(\bm z^{(k)})\|_2^2
    \\
    &\leq
    \left(
        \sum_{k=1}^{K'}
        \mathcal K_{{\rm D},M}(\bm z^{(k)})
    \right)
    \|Du\|_{L^2(\rho';\mathbb R^d)}^2.
\end{align*}
Because
\[
    \|u\|_{\rm EF}^2
    \geq
    w_{\rm F}\|Du\|_{L^2(\rho';\mathbb R^d)}^2,
\]
it follows that
\[
    \frac{
        w_{\rm F}
        \|\mathcal O_{\rm F}(u;\mathcal Z')\|_2^2
    }{
        \|u\|_{\rm EF}^2
    }
    \leq
    K'C_{\rm D}(U_M).
\]
Taking the suprema then proves
\[
\Lambda_{\rm F}^{\rm EF}(U_M)
    \leq
    \frac{K'}{w_{\rm F}}C_{\rm D}(U_M).
\]
\end{proof}

\subsection{Expected-case (\cref{prop:expected_observation_size})}\label{app:expected_observation_size}
\begin{proof}
For $u=\sum_{m=1}^M\theta_m\varphi_m$, the population energy geometry
gives
\[
    \|u\|_{\rm E}^2=\|\theta\|_2^2,
    \qquad
    \mathcal O_{\rm E}(u;\mathcal Z)
    =
    \sum_{m=1}^M
    \theta_m\Phi_m^{\rm E}(\mathcal Z).
\]
Hence
\[
    \sup_{\substack{u\in U_M\\\|u\|_{\rm E}=1}}
    \left|
        \mathcal O_{\rm E}(u;\mathcal Z)
    \right|^2
    =
    \sum_{m=1}^M
    |\Phi_m^{\rm E}(\mathcal Z)|^2.
\]
Taking expectations and using the orthonormality of the averaged energy
features gives
\[
    \mathbb E_{\mathcal Z}
    \left[
        \sup_{\substack{u\in U_M\\\|u\|_{\rm E}=1}}
        \left|
            \mathcal O_{\rm E}(u;\mathcal Z)
        \right|^2
    \right]
    =
    \sum_{m=1}^M
    \mathbb E_{\mathcal Z}
    |\Phi_m^{\rm E}(\mathcal Z)|^2
    =
    M.
\]
For the force observations, define
\[
    D(U_M)
    :=
    \{Du:u\in U_M\}
    =
    \operatorname{span}
    \{D\varphi_1,\ldots,D\varphi_M\}
    \subset L^2(\rho';\mathbb R^d),
\]
and let
\[
    r_D:=\dim D(U_M).
\]
Since $\varphi_1$ is constant, $D\varphi_1=0$, and hence
$r_D\leq M-1$.

Let $\{\psi_j\}_{j=1}^{r_D}$ be an orthonormal basis of $D(U_M)$ in
$L^2(\rho';\mathbb R^d)$. For every $u\in U_M$, there exist coefficients
$a_1,\ldots,a_{r_D}$ such that
\[
    Du=\sum_{j=1}^{r_D}a_j\psi_j,
    \qquad
    \|Du\|_{L^2(\rho';\mathbb R^d)}^2
    =
    \sum_{j=1}^{r_D}|a_j|^2.
\]
If $\|u\|_{\rm EF}=1$, then
\[
    w_{\rm F}\sum_{j=1}^{r_D}|a_j|^2
    =
    w_{\rm F}\|Du\|_{L^2(\rho';\mathbb R^d)}^2
    \leq1.
\]

Write
$\mathcal Z'=\{\bm z^{(1)},\ldots,\bm z^{(K')}\}$. By linearity,
\begin{align*}
    \mathcal O_{\rm F}(u;\mathcal Z')
    &=
    \frac{1}{\sqrt{K'}}
    \sum_{k=1}^{K'}Du(\bm z^{(k)})
    \\
    &=
    \sum_{j=1}^{r_D}
    a_j
    \left(
        \frac{1}{\sqrt{K'}}
        \sum_{k=1}^{K'}\psi_j(\bm z^{(k)})
    \right).
\end{align*}
Using Cauchy-Schwarz, we obtain
\begin{align*}
    \left\|
        \mathcal O_{\rm F}(u;\mathcal Z')
    \right\|_2^2
    &\leq
\left(\sum_{j=1}^{r_D}|a_j|^2\right)
    \sum_{j=1}^{r_D}
    \left\|
        \frac{1}{\sqrt{K'}}
        \sum_{k=1}^{K'}
        \psi_j(\bm z^{(k)})
    \right\|_2^2
    \\
    &\leq \frac{1}{w_{\rm F}}
    \sum_{j=1}^{r_D}
    \left\|
        \frac{1}{\sqrt{K'}}
        \sum_{k=1}^{K'}
        \psi_j(\bm z^{(k)})
    \right\|_2^2.
\end{align*}
Since the RHS is independent of $u$, it follows that
\begin{align}
    \sup_{\substack{u\in U_M\\\|u\|_{\rm EF}=1}}
    \left\|
        \mathcal O_{\rm F}(u;\mathcal Z')
    \right\|_2^2
    \leq \frac{1}{w_{\rm F}}
    \sum_{j=1}^{r_D}
    \left\|
        \frac{1}{\sqrt{K'}}
        \sum_{k=1}^{K'}
        \psi_j(\bm z^{(k)})
    \right\|_2^2.
    \label{eq:average_force_intermediate}
\end{align}

Each $\psi_j$ is the derivative of a function in $U_M$. Inversion symmetry
therefore implies
\[
    \int\psi_j(\bm z)\,d\rho'(\bm z)=0.
\]
Using independence of the cluster coordinates, we obtain
\begin{align*}
    &\mathbb E_{\mathcal Z'\sim\rho_{K'}}
    \left\|
        \frac{1}{\sqrt{K'}}
        \sum_{k=1}^{K'}
        \psi_j(\bm z^{(k)})
    \right\|_2^2
    \\
    &\qquad=
    \frac{1}{K'}
    \sum_{k,\ell=1}^{K'}
    \mathbb E
    \left[
        \psi_j(\bm z^{(k)})
        \cdot
        \psi_j(\bm z^{(\ell)})
    \right]
    \\
    &\qquad=
    \frac{1}{K'}
    \sum_{k=1}^{K'}
    \mathbb E
    \left[
        \|\psi_j(\bm z^{(k)})\|_2^2
    \right]
    \\
    &\qquad=
    \|\psi_j\|_{L^2(\rho';\mathbb R^d)}^2
    =1.
\end{align*}
Here the terms with $k\neq\ell$ vanish by independence and the zero-mean
property of $\psi_j$.

Taking expectations in
\eqref{eq:average_force_intermediate} now gives
\begin{align*}
    \mathbb E_{\mathcal Z'\sim\rho_{K'}}
    \left[
        \sup_{\substack{u\in U_M\\\|u\|_{\rm EF}=1}}
        \left\|
            \mathcal O_{\rm F}(u;\mathcal Z')
        \right\|_2^2
    \right]
    &\leq \frac{1}{w_{\rm F}}
    \sum_{j=1}^{r_D}1
    \\
    &= \frac{r_D}{w_{\rm F}}
    \leq \frac{M-1}{w_{\rm F}}.
\end{align*}
\end{proof}

\subsection{Pair-wise Case (\cref{prop:expected_observation_size})}\label{app:pair_wise_variation}
\begin{proof}[Proof of Proposition~\ref{prop:legendre_variation}]
By \eqref{eq:energy_christoffel_features} and
\eqref{eq:legendre_basis},
\[
\begin{aligned}
    \mathcal K_{{\rm E},M}(\mathcal Z)
    &=
    \sum_{m=1}^M
    \left|\Phi_m^{\rm E}(\mathcal Z)\right|^2
    \\
    &=
    1+
    \sum_{n=1}^{M-1}
    \left|
        \frac{1}{\sqrt K}
        \sum_{k=1}^K
        \sqrt{2n+1}\,
        P_n\bigl(\xi(r^{(k)})\bigr)
    \right|^2.
\end{aligned}
\]
For each $n$, Cauchy-Schwarz gives
\[
    \left|
        \frac{1}{\sqrt K}
        \sum_{k=1}^K
        \sqrt{2n+1}\,
        P_n\bigl(\xi(r^{(k)})\bigr)
    \right|^2
    \leq
    \sum_{k=1}^K
        (2n+1)
        \left|
            P_n\bigl(\xi(r^{(k)})\bigr)
        \right|^2.
\]
Therefore,
\[
    \mathcal K_{{\rm E},M}(\mathcal Z)
    \leq
    1+
    \sum_{k=1}^K
    \sum_{n=1}^{M-1}
        (2n+1)
        \left|
            P_n\bigl(\xi(r^{(k)})\bigr)
        \right|^2.
\]
Since $|P_n(x)|\leq1$ for $x\in[-1,1]$,
\[
    \sum_{n=1}^{M-1}
        (2n+1)|P_n(x)|^2
    \leq
    \sum_{n=1}^{M-1}(2n+1)
    =
    M^2-1.
\]
Hence, for every cluster collection $\mathcal Z$,
\[
    \mathcal K_{{\rm E},M}(\mathcal Z)
    \leq
    1+K(M^2-1).
\]
Taking the essential supremum over $\mathcal Z$ yields
\[
    \Lambda_{\rm E}(U_M)
    \leq
    1+K(M^2-1).
\]
For the energy component of the combined observations,~\cref{prop:generic_variation} gives
\[
    \Lambda_{\rm E}^{\rm EF}(U_M)
    \leq
    \frac{1}{w_{\rm E}}
    \left[1+K(M^2-1)\right].
\]

It remains to bound the force variation. Since $u$ is radial,
\[
    Du(\bm z)=u'(r)\widehat{\bm z},
    \qquad
    \|Du(\bm z)\|_2=|u'(r)|.
\]
Every $u\in U_M$ has the form
\[
    u(r)=p(\xi(r)),
\]
up to its constant term, where $p$ is a polynomial of degree at most
$M-1$. Therefore,
\[
    u'(r)=\xi'(r)p'(\xi(r)),
\]
and $p'$ has degree at most $M-2$.

We first bound $p'$ pointwise. Expand it in the normalized Legendre basis:
\[
    p'(x)
    =
    \sum_{j=0}^{M-2}
    b_j\sqrt{2j+1}\,P_j(x).
\]
By Cauchy-Schwarz and $|P_j(x)|\leq1$,
\begin{align*}
    |p'(x)|^2
    &\leq
    \left(\sum_{j=0}^{M-2}|b_j|^2\right)
    \left(\sum_{j=0}^{M-2}(2j+1)|P_j(x)|^2\right)
    \\
    &\leq
    (M-1)^2
    \sum_{j=0}^{M-2}|b_j|^2.
\end{align*}
By orthonormality of the normalized Legendre polynomials,
\[
    \sum_{j=0}^{M-2}|b_j|^2
    =
    \int_{-1}^1|p'(x)|^2\,\frac{dx}{2}.
\]
Therefore,
\begin{equation}\label{eq:polynomial_pointwise_bound}
    |p'(x)|^2
    \leq
    (M-1)^2
    \int_{-1}^1|p'(t)|^2\,\frac{dt}{2}.
\end{equation}

Using $|\xi'(r)|\leq C_\xi$ and
\eqref{eq:polynomial_pointwise_bound}, we obtain
\[
    \|Du(\bm z)\|_2^2 = |u'(r)|^2
    \leq
    C_\xi^2(M-1)^2
    \int_{-1}^1|p'(x)|^2\,\frac{dx}{2}.
\]
On the other hand, we have
\[
    \|Du\|_{L^2(\rho';\mathbb R^d)}^2
    =
    \int_{r_{\min}}^{r_{\rm cut}}
        |\xi'(r)|^2
        |p'(\xi(r))|^2
    \,d\rho_r(r).
\]
Since $|\xi'(r)|\geq c_\xi$,
\[
    \|Du\|_{L^2(\rho';\mathbb R^d)}^2
    \geq
    c_\xi^2
    \int_{r_{\min}}^{r_{\rm cut}}
        |p'(\xi(r))|^2
    \,d\rho_r(r).
\]
Finally, because the pushforward of $\rho_r$ under $\xi$ is uniform,
\[
    \|Du\|_{L^2(\rho';\mathbb R^d)}^2
    \geq
    c_\xi^2
    \int_{-1}^1
        |p'(x)|^2\,\frac{dx}{2}.
\]
Combining these inequalities gives
\begin{equation}\label{eq:pair_derivative_bound}
    \|Du(\bm z)\|_2^2
    \leq
    \left(\frac{C_\xi}{c_\xi}\right)^2
    (M-1)^2
    \|Du\|_{L^2(\rho';\mathbb R^d)}^2.
\end{equation}

Now write
$\mathcal Z'=\{\bm z^{(1)},\ldots,\bm z^{(K')}\}$. By
Cauchy-Schwarz and \eqref{eq:pair_derivative_bound},
\begin{align*}
    \left\|
        \mathcal O_{\rm F}(u;\mathcal Z')
    \right\|_2^2
    &=
    \frac{1}{K'}
    \left\|
        \sum_{k=1}^{K'}Du(\bm z^{(k)})
    \right\|_2^2
    \\
    &\leq
    \sum_{k=1}^{K'}
        \|Du(\bm z^{(k)})\|_2^2
    \\
    &\leq
    K'
    \left(\frac{C_\xi}{c_\xi}\right)^2
    (M-1)^2
    \|Du\|_{L^2(\rho';\mathbb R^d)}^2.
\end{align*}
Finally,
\[
    \|u\|_{\rm EF}^2
    \geq
    w_{\rm F}\|Du\|_{L^2(\rho';\mathbb R^d)}^2.
\]
Thus, since $\|u\|_{\rm EF}=1$,
\[
    \left\|
        \mathcal O_{\rm F}(u;\mathcal Z')
    \right\|_2^2
    \leq
    \frac{K'}{w_{\rm F}}
    \left(\frac{C_\xi}{c_\xi}\right)^2
    (M-1)^2.
\]
Taking the suprema proves
\eqref{eq:legendre_force_variation}.
\end{proof}

\section{Atomic Cluster Expansion}
\label{sec:ace_appndx}

We construct the higher body-order approximation spaces using the atomic
cluster expansion (ACE)~\cite{Drautz2019,DUSSON2022110946}. For a
cluster
\[
    \bm z=(\bm z_1,\ldots,\bm z_\nu),
    \qquad
    \bm z_j=r_j\widehat{\bm z}_j,
\]
the one-particle basis consists of products of radial functions and real
spherical harmonics,
\[
    \psi_{nlm}(\bm z_j)
    :=
    R_n(r_j)Y_l^m(\widehat{\bm z}_j).
\]
Summing over the neighbours gives the permutation-invariant features
\[
    A_{nlm}(\bm z)
    :=
    \sum_{j=1}^{\nu}\psi_{nlm}(\bm z_j).
\]
Products of these features are coupled to total angular momentum zero:
\[
    B_{\bm n\bm l q}(\bm z)
    :=
    \sum_{\bm m}
    C^{\bm l,q}_{\bm m}
    \prod_{\alpha=1}^{\nu}
    A_{n_\alpha l_\alpha m_\alpha}(\bm z),
\]
where \(C^{\bm l,q}_{\bm m}\) denotes the corresponding angular coupling
coefficients. Restricting to the scalar components gives a $O(3)$ invariant basis, while the neighbour sums give permutation
invariance.

In the numerical experiments, we use Legendre polynomials for the radial
basis and retain basis functions satisfying
\[
    \sum_{\alpha=1}^{\nu}(n_\alpha+l_\alpha)
    \leq D_{\rm tot}.
\]
Correlation orders \(\nu=2\) and \(\nu=3\) give the three- and four-body
spaces, respectively. The basis is generated using
\texttt{Polynomials4ML}\footnote{\url{https://github.com/ACEsuit/Polynomials4ML.jl}},
\texttt{EquivariantTensors}\footnote{\url{https://github.com/ACEsuit/EquivariantTensors.jl}.
}, and \texttt{SpheriCart}~\cite{sphericart}.

\bibliographystyle{siamplain}
\bibliography{references}
\end{document}